\documentclass[12pt, a4paper]{article}

\usepackage[utf8]{inputenc}
\usepackage{amsmath}   
\usepackage{amssymb}  
\usepackage{amsthm}   
\usepackage{mathtools} 

\usepackage{geometry}   
\usepackage{graphicx}   
\usepackage{subcaption}
\usepackage{natbib}     
\usepackage{booktabs}   
\usepackage{float}      

\usepackage{hyperref}   
\usepackage{cleveref}   

\newtheorem{theorem}{Theorem}[section]

\newtheorem{proposition}[theorem]{Proposition}
\newtheorem{corollary}[theorem]{Corollary}

\theoremstyle{definition}

\newtheorem{assumption}{Assumption}

\theoremstyle{remark}

\title{Contextual stochastic optimization using learned feature representations}

\author{%
\begin{minipage}{\textwidth}
\centering
{\large Ricardo Otero}\\[-0.25em]
{\footnotesize Research Center for Operations Management, Faculty of Economics and Business, KU Leuven}\\[-0.2em]
{\footnotesize Department of Industrial Engineering, Faculty of Engineering, Pontificia Universidad Javeriana}\\[-0.2em]
{\footnotesize \texttt{ricardo.otero@kuleuven.be}}\\[0.3em]
{\large Erik Demeulemeester}\\[-0.25em]
{\footnotesize Research Center for Operations Management, Faculty of Economics and Business, KU Leuven}\\[-0.2em]
{\footnotesize \texttt{erik.demeulemeester@kuleuven.be}}
\end{minipage}%
}
\date{}
\begin{document}

\maketitle

\begin{abstract}
Contextual stochastic optimization often uses weighted sample average approximation (SAA), assigning feature-dependent weights to historical scenarios. We examine the construction of these weights through the lens of learned feature representations, which map contextual feature vectors into spaces or structures where contexts can be compared. We develop representation-weighted contextual SAA, an estimate-then-optimize framework that separates feature representations, dissimilarities, and weight constructions. This separation makes the feature representation an explicit object of methodological study and broadens the range of learning methods available for contextual prescriptions to include methods that do not themselves construct scenario weights. The framework also allows representation-based weights to be combined with residual-based SAA. We derive finite-sample bounds on how closely the representation-weighted objective approximates the true conditional expected cost, accounting for errors in estimating the dissimilarities. For kernel weight constructions, we establish conditions for consistent prescriptions. Synthetic experiments show that learning how to represent and compare contexts can substantially improve the prescription quality relative to conventional contextual weighting methods, particularly in nonlinear settings. They also illustrate the benefits of exploring alternative ways of exploiting the representations and structures of existing fitted models. The experiments further show that informative learned representations can provide strong prescriptions without residual-based adjustments. A contextual electricity dispatch application using real-world data shows that representation-based weights can further improve residual-based SAA, even when the residual adjustments alone are highly effective. The framework thus provides a basis for exploiting advances in representation learning and investigating how representations and dissimilarities can be learned for contextual prescriptions.
\end{abstract}


\section{Introduction}

We study decision-making under uncertainty in settings where part of the variability in the uncertain model parameters can be explained by the features observed before the decision is made. An inventory planner, for example, may use weather forecasts, planned promotions, and recent sales before choosing stocking levels, as these features provide information about the demand in the upcoming planning period. Contextual stochastic optimization provides a framework for such situations by adapting decisions to the available features.

This setting can be formalized by letting \(X\in\mathcal X\subseteq\mathbb R^{d_x}\) denote the vector of features observed before the decision and \(\xi\in\Xi\subseteq\mathbb R^{d_\xi}\) denote the uncertain model parameters. After observing \(X=x\), the decision maker chooses \(z\in\mathcal Z\subseteq\mathbb R^{d_z}\) to minimize the conditional expected cost
\begin{equation}
\label{eq:conditional}
v^\star(x)=\min_{z\in\mathcal Z}\mathbb E[c(z,\xi)\mid X=x],
\end{equation}
where \(c(z,\xi)\) is the cost of the decision \(z\) under the realization \(\xi\). The conditional distribution of \(\xi\) given the observed features is unknown, and only a historical sample \(\mathcal D_n=\{(x_i,\xi_i)\}_{i=1}^{n}\) is available to approximate this problem. When \(c(z,\xi)\) is affine in \(\xi\), the conditional mean \(\mathbb E[\xi\mid X=x]\) fully determines the expected cost, whereas for more general cost functions the expected cost may depend on other aspects of the conditional distribution. This creates a distinction between settings that require a point prediction and those that require broader distributional estimation \citep{integrated_qi}. We focus on the latter.

A common data-driven strategy is estimate-then-optimize: a conditional distribution is estimated from the historical observations and then used in the stochastic optimization model \citep{surveycontextual_Sadana}. A prominent class of nonparametric prescriptions constructs this distribution by assigning probabilities to historical outcomes using the observed features. The outcomes \(\xi_i\) serve as scenarios, and their probabilities determine their contribution to the estimated expected cost. In an inventory problem, for example, each historical demand is a possible demand scenario. Given nonnegative weights \(w_i(x)\) that sum to one, the resulting weighted sample average approximation (SAA) is
\begin{equation}
\label{eq:weighted-saa}
\widehat z_n(x)\in\arg\min_{z\in\mathcal Z}
\sum_{i=1}^{n}w_i(x)c(z,\xi_i).
\end{equation}
This weighted empirical formulation underlies the framework for predictive prescriptions proposed by \citet{prescriptive_bertsimas}, which interprets the outputs of several machine-learning methods, including nearest neighbors, decision trees, and random forests, as weights over historical observations. The assigned weights determine the influence of each historical outcome, while the optimization problem determines the decision.

In this paper, we examine the construction of the weights in \eqref{eq:weighted-saa} through the lens of learned feature representations. A learned feature representation is a mapping estimated from data that expresses the original features in a space or structure in which observations can be compared. A learning method can supply such a representation without itself specifying the weights over the historical observations. We construct these weights by choosing how to compare observations in the learned representation and how to translate those comparisons into weights, for example, using a kernel. This separation makes the feature representation an explicit object of methodological study in contextual stochastic optimization. Representations learned for other purposes, including those obtained from pretrained models, can therefore be exploited for contextual prescriptions without changing the original role of the learning method. Beyond broadening the range of available learning methods, this perspective opens methodological questions concerning what information representations should preserve, how they should be learned, and how observations should be compared within them.

For example, dimensionality-reduction methods such as principal component analysis (PCA) and supervised nonlinear variants provide learned coordinates that can be used to construct the scenario weights \citep{dimensionality_ghojogh,spca_barshan}. Beyond these methods, our framework can accommodate neural-network embeddings extracted from hidden layers, allowing us to exploit expressive architectures capable of capturing complex relationships between the features and the uncertain parameters \citep{representationlearning_bengio}. The same fitted predictive model can thus be used to predict the uncertain parameters and, through its learned representations, to construct contextual scenario weights without modifying its original training procedure. 

Beyond expanding the range of learning methods, the representation perspective also allows alternative weight constructions within existing methods. Even when the representation is fixed, different dissimilarities can be used to compare observations. For vector representations, Euclidean \((\ell_2)\) and Manhattan \((\ell_1)\) distances provide different ways of comparing the same coordinates, or the dissimilarity itself can also be learned from data using metric learning \citep{mlkr_weinberger}. Our proposed perspective also reveals alternative constructions for fitted regression trees. Conventional tree weights use terminal leaves as representations, assigning equal weights to observations in the same leaf \citep{prescriptive_bertsimas}. However, fitted trees also provide richer representations through their routing paths, which record the nodes visited from the root to each terminal leaf. An alternative is to construct weights from these routing paths, using the number of edges connecting terminal leaves as the dissimilarity. This allows historical observations from other leaves to receive positive weights. Our experiments show that constructing the weights using learned metrics can improve the prescription quality relative to using distances in the original feature space, while weights based on routing paths can improve on conventional leaf weights for individual trees. The latter illustrates how alternative contextual prescriptions can be obtained from the same fitted model by changing how its learned structure is exploited.

A complementary way to incorporate the features is to adjust the historical outcomes before using them as scenarios. Residual-based approaches, for example, combine a prediction for the observed features with the prediction errors from historical observations, while retaining equal weights across the resulting scenarios \citep{residualtree_ban,predictivestochastic_deng,residualsaa_kannan}. Our framework naturally accommodates these transformations while also allowing representation-based weights to be assigned to the adjusted scenarios. This combination is particularly relevant when a transformation removes some contextual differences but leaves others that can be captured by learned representations and exploited through the resulting weights.

Our main contribution is the representation-weighted contextual SAA framework, which makes the representation of contextual information and the construction of dissimilarities explicit methodological choices in contextual stochastic optimization. We derive finite-sample bounds on how closely the representation-weighted objective approximates the true conditional expected cost, accounting for errors in estimating the dissimilarities. For kernel weight constructions, we establish conditions under which the resulting prescriptions are consistent. Drawing on results from the literature, we verify the dissimilarity estimation condition for PCA and its supervised and kernel variants, derive a perturbation bound for estimated Mahalanobis dissimilarities, and establish consistency for routing-path weights under explicit tree-growth conditions. Our synthetic experiments and a contextual electricity dispatch application using real-world data show that (i) learning how to represent and compare contexts can substantially improve prescription quality, particularly in nonlinear settings, with improvements also possible without reducing feature dimension, (ii) alternative prescriptions can be constructed from existing fitted models by exploiting their learned representations and structures, and (iii) representation-based weighting can complement residual-based SAA, even when the residual adjustments alone are highly effective. The framework thus provides a basis for developing representations specifically for contextual prescriptions and for exploiting representations learned for other purposes.
\subsection{Related work}

Our work belongs to the estimate-then-optimize literature in contextual optimization \citep{surveycontextual_Sadana} and focuses on risk-neutral contextual SAA with learned feature representations.

A central class of methods approximates the conditional expected cost by assigning feature-dependent weights to historical observations and using their outcomes as scenarios. Kernel and related local weighting methods have been used in data-driven dynamic programming \citep{robustdp_hanasusanto}, contextual newsvendor models \citep{bigdatanewsvendor_banrudin}, and other contextual prescriptions. \citet{prescriptive_bertsimas} place these ideas within a broader framework for predictive prescriptions by interpreting the outputs of several predictive methods as weights over historical observations. More recent work develops finite-sample and sample-complexity guarantees for Nadaraya--Watson contextual optimization \citep{regularizednw_wang}, characterizes context-dependent worst-case expected regret for weighted empirical-risk policies in contextual newsvendor problems \citep{contextualregret_besbes}, and studies statistical inference and computational-budget effects for weighted SAA \citep{weightedsaa_inference_wangzhang}. We build on this literature by studying the construction of scenario weights through the lens of learned feature representations. This perspective makes the underlying representations and dissimilarities explicit, allowing us to examine alternative constructions and their consequences for prescription quality.

Of particular relevance, \citet{regularizednw_wang} analyze a kernel weight construction after projecting the original features onto an estimated principal-component subspace to improve the sample complexity of contextual prescriptions. Their analysis accounts for the estimation of this projection and assumes that the retained components preserve the conditional information relevant to the uncertain parameters. Our work explores the broader role of learned feature representations and dissimilarities in contextual prescriptions. In our framework, dimensionality reduction is one possible use of learned representations rather than a requirement. We also examine how representation-based weighting interacts with scenario transformations.

A complementary literature provides methods for learning feature representations and dissimilarities. These methods differ in the information they seek to preserve and the relationships they can capture. Dimensionality-reduction methods learn lower-dimensional representations through linear or nonlinear transformations, either in an unsupervised manner, preserving structure in the contextual features, or in a supervised manner, using their relationship with the uncertain parameters \citep{dimensionality_ghojogh,isomap_tenenbaum,spca_barshan,kernelsdr_fukumizu}. Neural networks provide representations through the intermediate activations learned during model fitting \citep{representationlearning_bengio}, while metric-learning methods estimate dissimilarities using prediction, neighborhood, or other learning objectives \citep{metriclearning_bellet,mlkr_weinberger}. These methods therefore provide different ways to represent and compare observations. We investigate how their use in constructing scenario weights affects the contextual prescription quality.

Tree-based methods provide another way to learn representations through partitions of the feature space. Related work studies the statistical properties of data-dependent partitions and forest estimators \citep{partition_nobel,generalizedforest_athey,distributionalforest_cevid}. Conventional weights constructed from regression trees use terminal-leaf membership to determine which historical observations receive a positive weight \citep{prescriptive_bertsimas}. From a representation perspective, this is one possible use of the fitted tree structure. Our routing construction instead uses the relationships between the leaves to construct alternative weights without changing the fitted tree. We establish consistency for this construction under explicit tree-growth and localization conditions.

A separate stream of studies adapts the historical outcomes to the observed features. The residual-tree method of \citet{residualtree_ban} and predictive stochastic programming \citep{predictivestochastic_deng} construct scenarios by combining a prediction for the observed features with errors estimated from historical observations. \citet{residualsaa_kannan,residualdro_kannan} further develop residual-based SAA and distributionally robust formulations, including heteroscedastic location-scale constructions. These approaches accommodate different regression procedures for estimating the transformations. While residual-based SAA assigns equal weights to the adjusted scenarios, our framework lets their weights depend on the learned feature representations.

Another research direction incorporates the downstream optimization problem into model estimation. Task-based learning fits probabilistic models through the stochastic optimization objective \citep{taskbased_donti}, smart predict-then-optimize introduces a decision-aware surrogate loss \citep{smart_elmachtoub}, stochastic optimization forests select splits according to downstream decision quality \citep{stochasticforest_kallus}, integrated conditional estimation-optimization estimates conditional distributions using an objective aligned with the subsequent optimization problem \citep{integrated_qi}, and bilevel formulations directly account for downstream decisions during model training \citep{jimenez2026bilevel}. We instead focus on estimate-then-optimize, where feature representations are learned using statistical or prediction objectives, separately from the stochastic program. This separation allows us to exploit existing learning procedures and study the properties of their representations for contextual prescriptions.

The remainder of the paper is organized as follows. Section~\ref{sec:representation-weighted_SAA} presents representation-weighted contextual SAA and its main constructions. Section~\ref{sec:statistical_guarantees} develops the statistical analysis. Section~\ref{sec:computational_experiments} investigates the representation and scenario construction choices through controlled computational experiments, and presents the contextual electricity dispatch application. Finally, Section~\ref{sec:conclusions} concludes the paper by summarizing the main findings and their implications.

\section{Representation-weighted contextual SAA}
\label{sec:representation-weighted_SAA}
We now formalize the construction of contextual prescriptions from learned feature representations. The formulation distinguishes the representation used to express the features, the dissimilarity used to compare observations, and the mechanism used to construct the scenario weights. The historical outcomes may be used directly or transformed before entering the weighted SAA objective.

Let
\[
\phi:\mathcal X\rightarrow\mathcal H
\]
denote a feature representation, and let \(d\) be a nonnegative dissimilarity on the representation space \(\mathcal H\). The representation and the dissimilarity may be specified, estimated from the available data, or obtained from previously fitted models. For vector representations, \(\mathcal H\subseteq\mathbb R^{d_\phi}\), where \(d_\phi\) denotes the representation dimension which may be smaller than the original feature dimension \(d_x\). Together, the representation and the dissimilarity induce
\begin{equation}
\label{eq:dissimilarity}
\rho(x,x')
=
d\bigl(\phi(x),\phi(x')\bigr),
\qquad x,x'\in\mathcal X.
\end{equation}
The identity representation \(\phi(x)=x\) recovers the comparisons made in the original feature space.

In our main formulation, we construct the scenario weights using a nonnegative kernel \(K\) and a bandwidth \(h>0\). For an observed feature vector \(x\) and the historical observations \(\mathcal D_n=\{(x_i,\xi_i)\}_{i=1}^{n}\), observation \(i\) receives the weight
\begin{equation}
\label{eq:kernel-weights}
w_i(x)
=
\frac{K\bigl(\rho(x,x_i)/h\bigr)}
{\sum_{j=1}^{n}K\bigl(\rho(x,x_j)/h\bigr)},
\qquad i=1,\ldots,n.
\end{equation}
The kernel translates the scaled dissimilarities into weights through the normalization in \eqref{eq:kernel-weights}, while the bandwidth controls the scale at which the observations are compared. Provided that the denominator is positive, the weights are nonnegative and sum to one. Our theoretical analysis considers bounded kernels with compact support and the Gaussian kernel.

The historical outcome \(\xi_i\) may also be transformed using the observed features. Let \(\widehat T_x\) denote an optional scenario transformation and define
\begin{equation}
\label{eq:transformation}
\widetilde{\xi}_i(x)
=
\widehat T_x(x_i,\xi_i).
\end{equation}
For example, a fitted predictor \(\widehat m\) of the conditional mean gives the additive residual construction
\begin{equation}
\label{eq:residual}
\widetilde{\xi}_i(x)
=
\widehat m(x)+\xi_i-\widehat m(x_i).
\end{equation}
Here, the historical prediction error \(\xi_i-\widehat m(x_i)\) is added to the prediction at the observed feature vector \(x\). Applying equal weights to these adjusted outcomes recovers the additive residual-based SAA construction \citep{residualsaa_kannan}. When no transformation is used, \(\widehat T_x(x_i,\xi_i)=\xi_i\).

Combining the weights and the optional transformation gives the representation-weighted contextual SAA prescription
\begin{equation}
\label{eq:representation-saa}
\widehat z_n(x)
\in
\arg\min_{z\in\mathcal Z}
\sum_{i=1}^{n}
w_i(x)c\bigl(z,\widetilde{\xi}_i(x)\bigr).
\end{equation}
The weights determine the influence of each historical observation, while the transformation determines the outcome represented by its scenario. The representation and the transformation can be estimated separately, so a representation learned from the original observations can also be used to construct the weights for the transformed scenarios.

\begin{figure}[htbp]
    \centering
    \includegraphics[width=\linewidth]{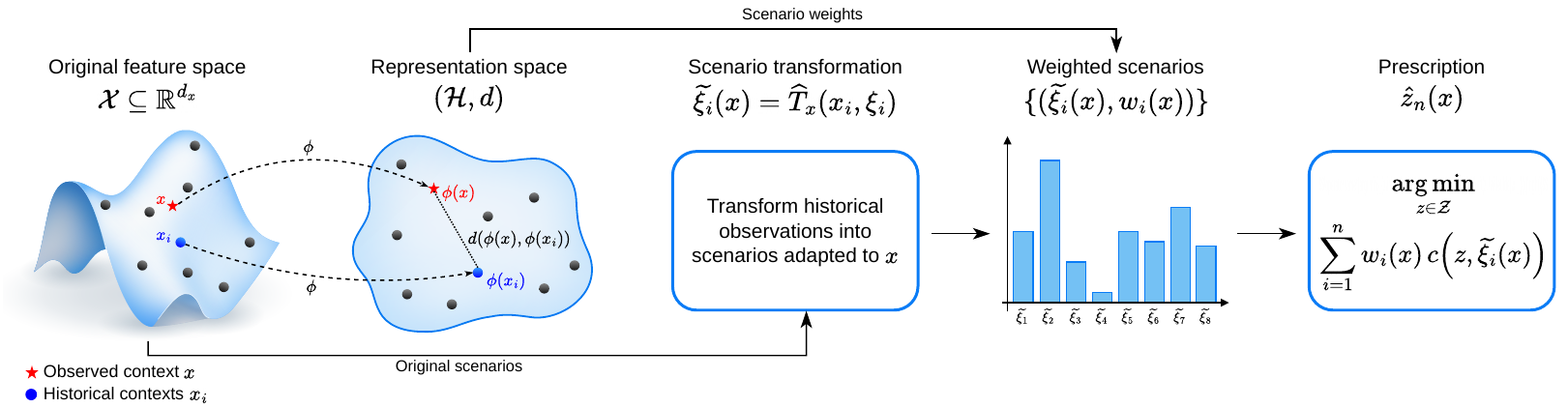}
    \caption{Representation-weighted contextual SAA.}
    \label{fig:framework}
\end{figure}

Figure~\ref{fig:framework} summarizes the components of the framework and their roles in the resulting prescription.

Table~\ref{tab:constructions} summarizes representative combinations of the feature representations, the dissimilarities, and the weight constructions. Alongside the kernel construction in \eqref{eq:kernel-weights}, it includes nearest-neighbor weighting, which assigns an equal weight \(1/k\) to the \(k\) historical observations closest to the observed feature vector under the induced dissimilarity \(\rho\), and zero weight to the others. The detailed definitions are provided in Appendix A.

\begin{table}[H]
\centering
\caption{Representative constructions within representation-weighted contextual SAA.}
\label{tab:constructions}
\footnotesize
\setlength{\tabcolsep}{4pt}
\renewcommand{\arraystretch}{1.2}
\begin{tabular}{
    p{\dimexpr0.18\linewidth-2\tabcolsep\relax}
    p{\dimexpr0.27\linewidth-2\tabcolsep\relax}
    p{\dimexpr0.27\linewidth-2\tabcolsep\relax}
    p{\dimexpr0.28\linewidth-2\tabcolsep\relax}
}
\hline
\textbf{Construction} & \textbf{Feature representation} & \textbf{Dissimilarity} & \textbf{Weight construction} \\
\hline
Original features & Identity, \(\phi(x)=x\) & Prescribed on \(\mathcal X\), such as the Euclidean distance & Kernel or nearest neighbors \\
Dimensionality reduction & Learned coordinates \(\phi(x)\) & Prescribed on \(\mathcal H\), such as the Euclidean distance & Kernel or nearest neighbors \\
Metric learning & Given representation \(\phi(x)\) & Learned on \(\mathcal H\) & Kernel or nearest neighbors \\
Tree: leaf & Terminal leaf, \(\phi(x)=L(x)\) & Membership dissimilarity & Uniform within the target leaf \\
Tree: routing & Routing-path representation \(\phi(x)\) & Routing-path distance & Kernel or nearest neighbors \\
Tree ensemble & Collection of tree-level representations & Tree-level dissimilarities & Mixture of separately normalized tree-level weights \\
Neural network & Hidden-layer representation \(\phi(x)\) & Prescribed on \(\mathcal H\), such as the Euclidean distance & Kernel or nearest neighbors \\
Grouping & Group label, \(\phi(x)=G(x)\) & Membership dissimilarity & Uniform within the target group \\
\hline
\end{tabular}
\end{table}

This flexibility allows the decision maker to construct alternative prescriptions by combining the available representations with different dissimilarities and weight constructions, using either the original historical outcomes or the transformed scenarios. Section~\ref{sec:statistical_guarantees} studies when the weighted objective in \eqref{eq:representation-saa} accurately approximates the true conditional expected cost and how its approximation error affects the quality of the resulting prescription.

\section{Statistical guarantees}
\label{sec:statistical_guarantees}

For the statistical analysis, fix the observed feature vector \(x\in\mathcal X\) and consider two independent i.i.d. samples from the joint distribution of \((X,\xi)\): a fitting sample \(\mathcal D_N^{\mathrm{fit}}\) and an SAA sample \(\mathcal D_n^{\mathrm{saa}}\). The fitting sample is used to estimate and tune the representation, dissimilarity, scenario transformation, and weight construction. The weights may additionally depend on the SAA features, but not on the SAA outcomes.

We denote the scenario weights by \(w_{N,n,i}(x)\) and the scenarios defined in Equation~\ref{eq:transformation} by \(\widetilde\xi_{N,i}(x)\), making explicit their dependence on the fitting sample size \(N\) and, for the weights, also on the SAA sample size \(n\). For \(z\in\mathcal Z\), define
\begin{equation}
\label{eq:conditional-objectives}
\widehat F_{N,n}(z\mid x)
=
\sum_{i=1}^{n}
w_{N,n,i}(x)c\bigl(z,\widetilde\xi_{N,i}(x)\bigr),
\qquad
F(z\mid x)
=
\mathbb E[c(z,\xi)\mid X=x],
\end{equation}
where the weights are nonnegative and sum to one. Our goal is to understand when the representation-weighted objective approximates the true conditional expected cost over \(\mathcal Z\) and therefore yields a good contextual prescription.

Standard SAA theory provides guarantees for this approximation when scenarios are sampled from \(\xi\mid X=x\). Our analysis additionally bounds the discrepancy introduced by using historical observations with representation-based weights instead of samples from this conditional distribution. We first derive finite-sample guarantees for general weight constructions satisfying the stated assumptions. We then establish consistency conditions specifically for kernel weight constructions and, finally, analyze the effect of scenario transformations.

\subsection{The role of the learned dissimilarity}
\label{sec:role-dissimilarity}

The representation-weighted objective approximates the unknown conditional expected cost at the observed features \(x\) using historical outcomes associated with other feature vectors \(x'\). We measure the discrepancy at each \(x'\) by taking the largest difference over the feasible decisions:
\[
\Gamma_x(x')
=
\sup_{z\in\mathcal Z}
\left|F(z\mid x')-F(z\mid x)\right|.
\]

To assess how these discrepancies affect the approximation, we consider their weighted average, \(\sum_{i=1}^{n}w_{N,n,i}(x)\Gamma_x(X_i)\). Since the weights favor observations with smaller dissimilarities, we relate these dissimilarities to the differences in the conditional expected costs.

Let \(\phi_0:\mathcal X\to\mathcal H_0\) denote the population feature representation and \(d_0\) a dissimilarity on \(\mathcal H_0\), inducing \(\rho_0(x,x')=d_0(\phi_0(x),\phi_0(x'))\). We assume that there exists a finite constant \(L_\rho\geq0\) such that
\begin{equation}
\label{eq:population_relevance_condition}
\Gamma_x(x')
\leq
L_\rho\rho_0(x,x'),
\qquad x'\in\mathcal X.
\end{equation}

Condition~\eqref{eq:population_relevance_condition} characterizes a sufficient property of the population representation and the dissimilarity for contextual prescriptions: contexts that are close under \(\rho_0\) must have similar conditional expected costs across all feasible decisions. Preserving the conditional mean may be insufficient for this property, while preserving the entire conditional distribution may be unnecessary for a given optimization problem. A distributional condition could nevertheless be useful when considering a broader class of cost functions. For example, if the 1-Wasserstein distance between the conditional distributions is bounded by a constant multiple of \(\rho_0\), then \eqref{eq:population_relevance_condition} holds for any cost function that is Lipschitz continuous in the uncertainty, uniformly over the feasible decisions. Related conditions appear in \citet[Assumption~(A3)]{kernelgen_wang} for conditional expected costs and in \citet{bertsimas_dynamic} for conditional distributions, both using distances in the original feature space.

Condition~\eqref{eq:population_relevance_condition} concerns the relevance of the population dissimilarity, independently of the accuracy of its estimation. In practice, the scenario weights are constructed using a representation and a dissimilarity estimated from the fitting sample, yielding the fitted counterpart \(\widehat{\rho}_N\) of \(\rho_0\). For each observation in the SAA sample, the population and fitted dissimilarities satisfy
\[
\rho_0(x,X_i)
\leq
\widehat\rho_N(x,X_i)+\Delta_{N,n}(x),
\qquad i=1,\ldots,n,
\]
where
\begin{equation}
\label{eq:fitted_dissimilarity_error}
\Delta_{N,n}(x)
=
\max_{1\leq i\leq n}
\left|
\widehat\rho_N(x,X_i)-\rho_0(x,X_i)
\right|
\end{equation}
measures the largest estimation error across these comparisons. The convergence behavior of \(\Delta_{N,n}(x)\) depends on the estimation method and the assumptions imposed on the learned representation and dissimilarity. We examine this in more detail for some representations in Appendix \ref{app:method_verification}.

Combining these bounds and taking the weighted sum yields
\[
\begin{aligned}
\sum_{i=1}^{n}w_{N,n,i}(x)\Gamma_x(X_i)
&\leq
L_\rho\left(
\sum_{i=1}^{n}w_{N,n,i}(x)\widehat\rho_N(x,X_i)
+\Delta_{N,n}(x)
\right)\\
&=
L_\rho\bigl(\widehat R_{N,n}(x)+\Delta_{N,n}(x)\bigr),
\end{aligned}
\]
where \(\widehat R_{N,n}(x)=\sum_{i=1}^{n}w_{N,n,i}(x)\widehat\rho_N(x,X_i)\) is the weighted radius under the fitted dissimilarity. This separates the proximity induced by the fitted dissimilarity from the error in estimating it.

Combining this result with a standard uniform concentration argument in stochastic programming \citep[Chapter~5]{saa_shapiro} gives a bound on the difference between the representation-weighted objective and the true conditional expected cost. We retain the independent fitting and SAA samples introduced above and assume that the feasible set is nonempty and compact, the cost function is Lipschitz continuous in the decision and the uncertain parameters, and the conditional fluctuations of the scenario costs satisfy a uniform sub-Gaussian tail condition. The complete assumptions and proof are provided in Appendix B.

\begin{theorem}[Bound using the fitted representation]
\label{cor:fitted_representation_bound}
Under the conditions described above, stated formally in Assumptions~\ref{ass:independent_samples}--\ref{ass:relevance_dissimilarity} of Appendix~\ref{app:proofs_statistical_guarantees}, for every \(\eta>0\) and \(\delta\in(0,1)\), with conditional probability at least \(1-\delta\) given the fitting sample and the SAA features,
\begin{equation}
\label{eq:fitted_representation_bound}
\begin{aligned}
\sup_{z\in\mathcal Z}
\left|
\widehat F_{N,n}(z\mid x)-F(z\mid x)
\right|
&\leq
E^T_{N,n}(x)
+
L_\rho\widehat R_{N,n}(x)
+
L_\rho\Delta_{N,n}(x)\\
&\quad+
2L_z\eta
+
\bar\tau_x
\sqrt{
2\log\!\left(
\frac{2\mathcal N_{\mathcal Z}(\eta)}{\delta}
\right)
\sum_{i=1}^{n}w_{N,n,i}(x)^2
}.
\end{aligned}
\end{equation}
\end{theorem}

Here, \(E^T_{N,n}(x)\) measures the error from estimating the scenario transformation, which is analyzed in Section~\ref{sec:scenario_transformations}. When no transformation is used, \(\widetilde\xi_{N,i}(x)=\xi_i\) and \(E^T_{N,n}(x)=0\). The constant \(L_z\) is the Lipschitz constant of the cost function with respect to the decision, \(\bar\tau_x^2\) is the conditional variance proxy in the sub-Gaussian assumption, and \(\mathcal N_{\mathcal Z}(\eta)\) denotes the \(\eta\)-covering number of \(\mathcal Z\).

The theorem shows how the weighted radius, dissimilarity estimation error, transformation error, and sampling variability contribute to the approximation error. Let \(\widehat z_{N,n}(x)\) minimize the estimated objective \(\widehat F_{N,n}(\cdot\mid x)\). By the optimality of \(\widehat z_{N,n}(x)\) and the triangle inequality,
\begin{equation}
\label{eq:lfrB-regret}
0\leq F\bigl(\widehat z_{N,n}(x)\mid x\bigr)-v^\star(x)\leq 2\sup_{z\in\mathcal Z}\left|\widehat F_{N,n}(z\mid x)-F(z\mid x)\right|.
\end{equation}
Thus, uniform objective consistency implies vanishing prescription regret. We next establish sufficient conditions for these properties under kernel weight constructions.

\subsection{Kernel weights in a fitted representation}
\label{sec:kernel_specialization}

\begin{figure}[t]
    \centering
    \includegraphics[width=0.32\linewidth]{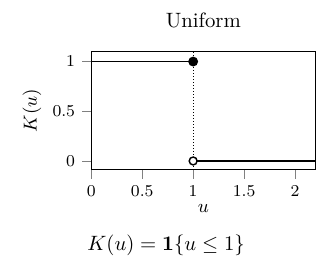}%
    \hfill
    \includegraphics[width=0.32\linewidth]{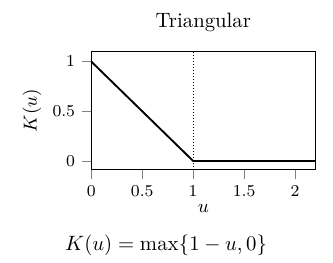}%
    \hfill
    \includegraphics[width=0.32\linewidth]{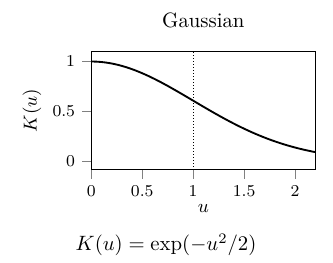}
    \caption{Kernel functions applied to the scaled dissimilarity
    \(u=\widehat\rho_N(x,X_i)/h_N\). The dotted lines mark \(u=1\),
    corresponding to a fitted dissimilarity of \(h_N\). Here, \(\mathbf{1}\{\cdot\}\) denotes the indicator function.}
    \label{fig:kernel_examples}
\end{figure}

In the kernel construction, the historical observations are compared with the observed features \(x\) through their fitted dissimilarities, which the kernel converts into weights. The kernel is evaluated at the scaled dissimilarity \(u=\widehat\rho_N(x,X_i)/h_N\), where the bandwidth \(h_N>0\) sets the scale of these comparisons. Our analysis below considers nonnegative bounded kernels with compact support, including the uniform and triangular kernels, as well as the Gaussian kernel. Figure~\ref{fig:kernel_examples} illustrates these examples.

A smaller bandwidth makes the approximation more local but can increase the sampling variability by concentrating the weights on fewer observations. Classical consistency arguments for kernel regression balance these effects by letting the bandwidth shrink as the sample grows, while retaining an increasing number of observations in shrinking neighborhoods; see, for example, \citet[Chapters~4 and~6]{functionaldata_ferratyvieu}.

We summarize how many scenarios effectively support the approximation using the effective sample size
\[
n_{\mathrm{eff},N,n}(x)
=
\frac{1}{\sum_{i=1}^{n}w_{N,n,i}(x)^2},
\]
which equals \(n\) under equal weights and decreases when a few scenarios carry most of the weight.

To ensure that enough SAA observations remain near \(x\) as the bandwidth shrinks, we use the local mass condition in Assumption~\ref{ass:local_mass} of Appendix~\ref{app:kernel_consistency}. It requires that, for some \(q>0\) and all sufficiently small \(s>0\), the probability of finding an observation within dissimilarity \(s\) of \(x\) is at least a positive constant times \(s^q\). Thus, for sufficiently small \(h_N\), the expected number of SAA observations within dissimilarity \(h_N\) of \(x\), conditional on the fitting sample, is at least a constant times \(nh_N^q\). 

For example, under Euclidean dissimilarity in a \(d_\phi\)-dimensional representation, the condition holds with \(q=d_\phi\) if \(\phi_0(X)\) has a density that is continuous and positive at \(\phi_0(x)\). When \(d_\phi\ll d_x\), the factor \(h_N^{d_\phi}\) decreases more slowly than \(h_N^{d_x}\) as the bandwidth shrinks. A representation that preserves the relevant information in fewer dimensions can therefore require fewer observations to maintain local support. Similar local mass conditions appear in kernel regression \citep[Chapter~13]{functionaldata_ferratyvieu} and contextual kernel analysis \citep[Assumption~(A2) and Proposition~1]{kernelgen_wang}.
The local mass condition characterizes proximity under the population dissimilarity \(\rho_0\), whereas the kernel weights are constructed using the fitted dissimilarity \(\widehat{\rho}_N\). Since the kernel compares dissimilarities relative to the bandwidth \(h_N\), we require the estimation error to vanish relative to this scale:
\[
\frac{\Delta_{N,n}(x)}{h_N}\xrightarrow{p}0,
\]
where \(\xrightarrow{p}\) denotes convergence in probability. 

Combining the conditions above with Theorem~\ref{cor:fitted_representation_bound} gives the following consistency result, establishing that the prescription regret converges to zero in probability.

\begin{corollary}[Consistency of kernel prescriptions]
\label{cor:kernel_consistency}
Under the assumptions of Theorem~\ref{cor:fitted_representation_bound} and the additional conditions described above, which are stated in detail in Appendix~\ref{app:kernel_consistency}, let the bandwidth \(h_N\) be selected using the fitting sample. If, as \(N,n\to\infty\),
\[
h_N\xrightarrow{p}0,\qquad nh_N^q\xrightarrow{p}\infty,\qquad \frac{\Delta_{N,n}(x)}{h_N}\xrightarrow{p}0,\qquad E^T_{N,n}(x)\xrightarrow{p}0,
\]
then
\[
\sup_{z\in\mathcal Z}\left|\widehat F_{N,n}(z\mid x)-F(z\mid x)\right|\xrightarrow{p}0,
\qquad
F\bigl(\widehat z_{N,n}(x)\mid x\bigr)-v^\star(x)\xrightarrow{p}0.
\]
\end{corollary}

The proof is provided in Appendix~\ref{app:kernel_consistency}.

Under the stated assumptions, these results establish the consistency of the prescriptions when (i) the bandwidth shrinks while enough historical scenarios contribute to the approximation, (ii) the dissimilarity estimation error is negligible relative to the bandwidth, and (iii) the transformation estimation error vanishes. In Appendix~\ref{app:spectral_verification}, we give conditions under which the dissimilarity estimation requirement holds for PCA and its supervised and kernel variants. In Appendix~\ref{app:mahalanobis_verification}, we relate the estimation error of the Mahalanobis dissimilarity to the estimation error of the Mahalanobis matrix. In Appendix~\ref{app:tree_consistency}, we establish consistency for the weights based on routing paths under explicit conditions on the tree construction. We evaluate the neural representations only empirically in Section~\ref{sec:computational_experiments}.

\subsection{Scenario transformations}
\label{sec:scenario_transformations}

A scenario transformation \(\widehat T_{N,x}\) adjusts the historical outcomes using the observed features \(x\) before they enter the weighted objective. It can remove some systematic differences associated with the features, while the representation-based weights can exploit the remaining differences relevant to the conditional expected costs.

Let \(T_x^0\) denote the population counterpart of \(\widehat T_{N,x}\). We measure the differences that remain after the population transformation by
\[
\Gamma_x^T(x')
=
\sup_{z\in\mathcal Z}
\left|
\mathbb E[c(z,T_x^0(x',\xi))\mid X=x']
-
F(z\mid x)
\right|.
\]
This compares the conditional expected costs of the transformed outcomes originating from \(x'\) with the true conditional expected costs at \(x\). Using \(\Gamma_x^T\) instead of \(\Gamma_x\) in condition~\eqref{eq:population_relevance_condition} means that the population dissimilarity only needs to control the differences that remain after the transformation.

Estimating this transformation introduces a separate error. Let \(E^T_{N,n}(x)\) denote the largest absolute change in the weighted objective over the feasible decisions when the fitted transformation replaces its population counterpart. Writing \(\widetilde\xi_i^0(x)=T_x^0(X_i,\xi_i)\), the Lipschitz continuity of the cost function in the uncertain parameters gives
\[
E^T_{N,n}(x)
\leq
L_\xi
\sum_{i=1}^{n}w_{N,n,i}(x)
\left\|
\widetilde\xi_{N,i}(x)-\widetilde\xi_i^0(x)
\right\|.
\]
If the fitted transformation converges uniformly to its population counterpart, then \(E^T_{N,n}(x)\xrightarrow{p}0\), satisfying the transformation estimation condition in Corollary~\ref{cor:kernel_consistency}. For residual constructions, \citet{residualsaa_kannan} collect and present consistency results for several statistical and machine-learning regression methods under suitable assumptions.

For the additive residual construction in \eqref{eq:residual}, recentering subtracts the predicted conditional mean at the historical features and adds the prediction at \(x\). When the conditional distribution differs across contexts only through its mean, recentering with the true conditional mean recovers the conditional distribution at \(x\), leaving no differences for the representation-based weights to address. However, a misspecified mean predictor can leave differences after recentering, and even using the true conditional mean may leave differences in dispersion, shape, or tails that affect the conditional expected costs. A representation learned from the original observations may capture these remaining differences, allowing its scenario weights to complement recentering. We examine this in the computational experiments in Section~\ref{sec:computational_experiments}.

More general scenario transformations, such as location--scale constructions, can remove additional contextual differences by adjusting the conditional dispersion. However, as discussed by \citet{residualsaa_kannan}, estimating the conditional dispersion from the available data can be challenging and typically requires assumptions about its functional form. In their scale-adjusted experiments, they assume that this form is known. Representation-based weighting provides a complementary way to account for the differences that remain after a chosen transformation. We therefore focus our experiments on mean recentering and investigate whether the representation-based weights can exploit the remaining differences without explicitly modeling the conditional dispersion.

\section{Computational experiments}
\label{sec:computational_experiments}

We evaluate representation-weighted contextual SAA through controlled synthetic experiments and a contextual electricity dispatch application. The former allows us to examine how different representation and scenario constructions affect the prescription quality, while the latter evaluates their performance using real-world data.

\subsection{Synthetic experiments}
\label{sec:synthetic_experiments}
\label{sec:computational_common}

We use a fixed optimization problem while varying the relationships between the contextual features and demand, the learning methods, and the available data.

\paragraph{Problem description.}
We consider a two-stage capacity-expansion problem with four machines and four products. After observing the contextual features \(x\), a manufacturer chooses the available machine capacities \(t\), paying for any expansion beyond the initial capacities \(t^0\). The problem is formulated as
\begin{equation}
\label{eq:capacity_expansion_problem}
\min_{t^0\leq t\leq t^{\max}}
\left\{
\sum_{m\in\mathcal M}c_m^E(t_m-t_m^0)
+
\mathbb E[Q(t,\xi)\mid X=x]
\right\},
\end{equation}
where \(c_m^E\) is the expansion cost for machine \(m\). Once demand \(\xi\) is realized, the recourse problem allocates production across the machines and penalizes the unmet demand:
\begin{equation}
\label{eq:capacity_expansion_recourse}
\begin{aligned}
Q(t,\xi)=\min_{u,s\geq0}\quad
&\sum_{(j,m)\in\mathcal A}c_{jm}^P u_{jm}
+\sum_{j\in\mathcal J}c_j^S s_j\\
\text{s.t.}\quad
&\sum_{m:(j,m)\in\mathcal A}\gamma_{jm}u_{jm}+s_j
\geq\xi_j,
&&j\in\mathcal J,\\
&\sum_{j:(j,m)\in\mathcal A}u_{jm}
\leq t_m,
&&m\in\mathcal M.
\end{aligned}
\end{equation}

Here, \(u_{jm}\) is the operating time allocated to product \(j\) on machine \(m\), \(\gamma_{jm}\) is the number of units of product \(j\) produced per unit of operating time on machine \(m\), and \(s_j\) represents the unmet demand. The coefficients \(c^P\) and \(c^S\) denote production and shortage costs, respectively, and \(\mathcal A\) contains the feasible product-machine assignments. Each product can be produced on a primary or secondary machine, with secondary production being slower and more expensive. The numerical parameterization is provided in Appendix~\ref{app:synthetic_dgps}.

\paragraph{Data-generating processes.}
We consider four data-generating processes (DGPs) based on four Gaussian demand factors \(U\sim\mathcal N(0,I_4)\), with each component clipped to the interval\([-6,6]\). The demand for each product \(j\in\mathcal J\) is generated as
\begin{equation}
\label{eq:synthetic_demand_generation}
\xi_j
=
\left[
24\left(1+\alpha\widetilde g_j(U)\right)
+\eta_j(U)
\right]_+,
\qquad \alpha=0.18,
\end{equation}

where \([a]_+=\max\{a,0\}\) denotes the positive part, ensuring nonnegative demand, \(\widetilde g_j(U)\) is the standardized demand signal, and \(\eta_j(U)\) is a Gaussian noise term with possibly context-dependent variance. The demand signal is linear in the factors in the first two DGPs and nonlinear in the remaining two.

In the linear DGP, the demand factors are directly observed as the contextual features. 
In the latent DGPs, the factors are not directly observed but are mixed with independent Gaussian noise variables \(V\sim\mathcal N(0,I_{d_v})\), also clipped to \([-6,6]\), through an orthogonal transformation. Here, \(d_v=d_x-4\) denotes the number of noise dimensions, which are not a priori distinguishable from the informative dimensions. The informative directions have larger variance than the noise directions in the linear latent setting, whereas they have equal variance in the nonlinear latent settings. The final DGP additionally introduces a context-dependent noise variance. Table~\ref{tab:synthetic_dgps} summarizes these settings, with the complete data-generation procedures provided in Appendix~\ref{app:synthetic_dgps}.

\begin{table}[t]
\centering
\caption{Core data-generating processes and baseline feature dimensions.}
\label{tab:synthetic_dgps}
\footnotesize
\setlength{\tabcolsep}{2pt}
\renewcommand{\arraystretch}{1.15}
\begin{tabular}{@{}p{0.40\linewidth}rrll@{}}
\toprule
DGP & \(d_x\) & \(d_v\) & Signal & Noise \\
\midrule
Linear & 4 & 0 & Linear & Constant \\
Linear latent & 30 & 26 & Linear & Constant \\
Nonlinear latent & 16 & 12 & Nonlinear & Constant \\
Nonlinear latent, conditional variance
& 16 & 12 & Nonlinear & Variable \\
\bottomrule
\end{tabular}
\end{table}

\paragraph{Methods and benchmarks.}
We compare contextual prescriptions that differ in their feature representations, dissimilarities, and weight constructions. The \emph{Raw} method uses Euclidean distances between the standardized original features, corresponding to conventional contextual kernel weighting without a learned representation \citep{prescriptive_bertsimas, robustdp_hanasusanto}. \emph{PCA} projects these features onto their leading principal components and constructs kernel weights in the reduced space, following the dimensionality reduction approach of \citet{regularizednw_wang}. We also evaluate machine-learning methods that can be incorporated into contextual optimization through their learned feature representations. \emph{KSPCA} uses kernel supervised PCA with a polynomial kernel to learn nonlinear representations informed by demand \citep{spca_barshan}, providing a supervised nonlinear alternative to PCA. This representation kernel is distinct from the kernel used to construct scenario weights. To examine whether learning the dissimilarity between observations improves the prescriptions, \emph{MLKR} estimates a distance by minimizing the leave-one-out kernel regression prediction error \citep{mlkr_weinberger}. The \emph{Neural} method uses the final hidden-layer representation of a feed-forward network trained to predict demand, providing a more flexible nonlinear representation. For these methods, Euclidean distances in the representation space are converted into scenario weights using a kernel.

We also consider tree constructions. \emph{DT leaf} assigns equal weights to observations in the same terminal leaf, corresponding to conventional tree weighting \citep{prescriptive_bertsimas}. \emph{DT route} instead uses a kernel applied to the routing-path distance, allowing observations from different leaves to receive positive weights. Similarly, \emph{RF leaf} and \emph{RF route} average the corresponding normalized weights across trees. These comparisons use the same fitted trees or forests, changing only the weight construction to isolate its effect on the prescriptions.

We use the Gaussian kernel throughout the remaining experiments and examine the effect of alternative kernels in Appendix~\ref{app:kernel_comparison}.

We include two oracle benchmarks. The \emph{oracle prescription} constructs prescriptions using independent scenarios drawn from the true conditional distribution \(\xi\mid X=x\). The \emph{oracle representation} instead uses the true demand factors \(U\) in the linear DGPs and their nonlinear transformations \(\psi(U)\) in the nonlinear DGPs, but still relies on the Gaussian-kernel weights over historical scenarios to approximate the conditional distribution. We also include \emph{SAA}, which assigns equal weights to all historical scenarios.

\paragraph{Experimental protocol.}
Unless stated otherwise, we use \(N=1500\) fitting observations (1000 for training and 500 for validation) and an independent SAA sample of \(n=1000\) observations. We evaluate the prescriptions at 100 test feature vectors over 25 independent replications.

For each test feature vector \(x\), we fix the prescribed first-stage decision and estimate its conditional expected cost by solving the recourse problem for 500 independent demand realizations drawn from \(\xi\mid X=x\). The oracle prescription (OP) is approximated using 1000 additional conditional demand realizations. Let \(C_m^{(i)}\) denote the mean evaluated cost across the test feature vectors for method \(m\) in replication \(i\). We measure performance using the normalized gap
\begin{equation}
\label{eq:normalized_gap}
\mathrm{gap}_m^{(i)}
=
100\,
\frac{C_m^{(i)}-C_{\mathrm{OP}}^{(i)}}
{C_{\mathrm{SAA}}^{(i)}-C_{\mathrm{OP}}^{(i)}}.
\end{equation}

The oracle prescription corresponds to \(0\%\) and SAA to \(100\%\), with lower values indicating better performance. Below, the boxplots summarize the 25 replication-level normalized gaps, while the line plots show their mean with shaded pointwise 95\% confidence intervals across replications.

We use the validation observations to select the tuning parameters. The neural learning rate and weight decay and the tree minimum leaf size are selected using the prediction RMSE, while the representation dimension and Gaussian-kernel bandwidth are selected jointly using the energy score, which evaluates the estimated conditional distribution \citep{energyscore_gneiting}. The complete tuning procedures and implementation details are provided in Appendix~\ref{app:synt_training_details}.

\subsection{Effect of representation choice}
\label{sec:exp2_representation}

\begin{figure}[t]
    \centering
    \includegraphics[width=0.8\textwidth]{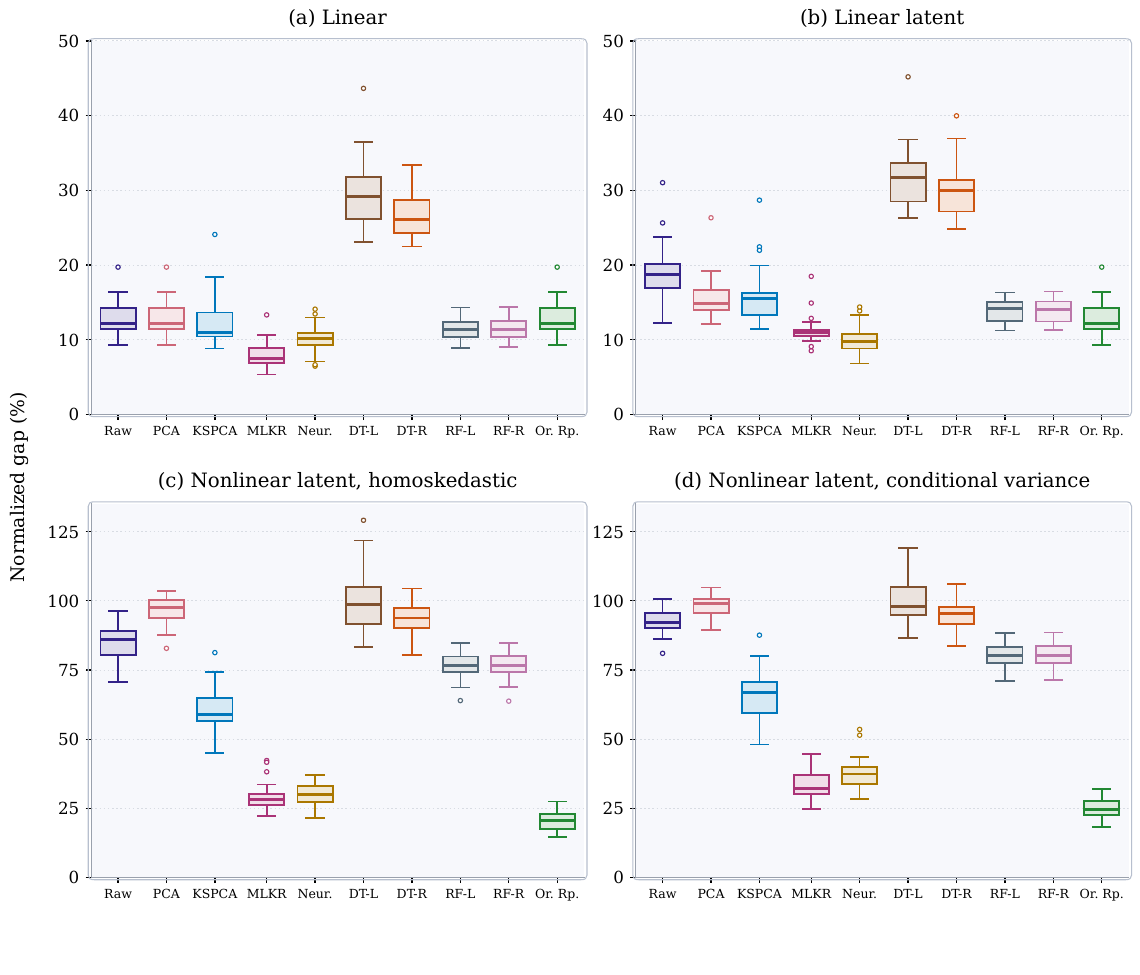}
    \caption{Normalized gap across the four DGPs.}
    \label{fig:exp2a}
\end{figure}

We first compare the prescription performance of different feature representations and weight constructions across the four DGPs. Figure~\ref{fig:exp2a} shows that the relative performance of the methods changes substantially with the demand structure. \emph{Raw} is competitive when the demand factors are directly observed, while \emph{PCA} improves relative to \emph{Raw} in the linear-latent setting, where the informative directions contribute more strongly to the variation in \(X\). This pattern reverses in the nonlinear latent settings, where \emph{PCA} performs close to \emph{SAA} and worse than \emph{Raw}. The additional representation methods considered in our framework provide substantial improvements, with \emph{KSPCA} offering intermediate gains and \emph{MLKR} and \emph{Neural} achieving particularly strong performance in the nonlinear latent settings. These results illustrate the benefit of incorporating a broader range of learning methods into contextual SAA. The gaps of the strongest methods are larger in the conditional-variance DGP, but their relative performance changes little. Note that the oracle representation can be outperformed, since knowing the true demand factors does not necessarily mean that Euclidean distances provide the best scenario weights.

The results also illustrate the value of alternative constructions within existing methods. \emph{MLKR} improves on \emph{Raw} across all four DGPs, highlighting the benefit of learning the distance rather than using Euclidean comparisons in the original feature space. Similarly, \emph{DT route} improves on \emph{DT leaf} across the four DGPs. This improvement is particularly relevant when a single tree is preferred for its interpretable splitting structure, since routing-path weighting improves prescriptions without changing the fitted tree. In contrast, \emph{RF leaf} and \emph{RF route} perform similarly. One possible explanation is that averaging weights across many trees already combines observations from different leaves, reducing the additional benefit of routing-path distances. Both forest constructions are competitive in the linear settings but deteriorate substantially in the nonlinear latent DGPs.

A further comparison fixing the representation dimension at \(d_\phi=4\) confirms that performance differences remain across learning methods, indicating that the choice of method matters beyond representation dimension alone. This comparison is reported for the two nonlinear DGPs in Appendix~\ref{app:common_dimension}.

\subsection{Effect of scenario transformations}
\label{sec:exp_transf}

We next consider how the additive residual construction interacts with representation-based weighting. We compare the original scenarios (\emph{Obs}) with scenarios recentered using either a LightGBM estimate of the conditional mean (\emph{E-Rec}) or a Monte Carlo approximation of the true conditional mean based on 10,000 conditional draws per feature vector (\emph{O-Rec}). Combining \emph{E-Rec} with equal weights recovers the residual-based SAA construction of \citet{residualsaa_kannan}. Alongside the two nonlinear latent DGPs, we consider variants with context-dependent positive shocks and a bimodal structure, described in Appendix~\ref{app:synthetic_dgps}. We report results for \emph{SAA}, \emph{Raw}, \emph{KSPCA}, \emph{MLKR}, \emph{Neural}, \emph{RF leaf}, and the oracle representation.

\begin{figure}[t]
    \centering
    \includegraphics[width=0.80\textwidth]{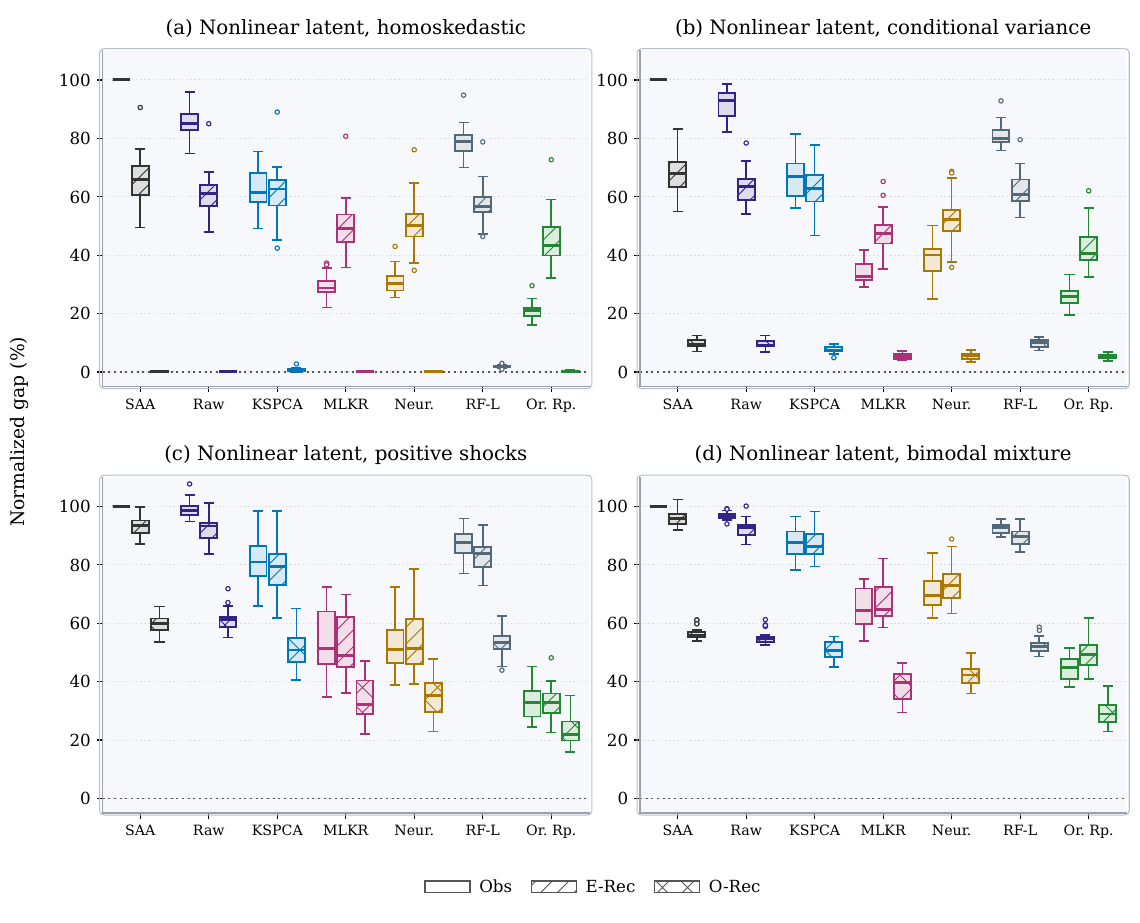}
    \caption{Normalized gap under original, estimated-mean-recentered, and oracle-mean-recentered scenarios across the nonlinear latent DGPs.}
    \label{fig:exp3}
\end{figure}

Figure~\ref{fig:exp3} shows that representation-based weights can improve on residual-based SAA with equal weights, while the benefit of recentering depends on the DGP and the learning method. \emph{E-Rec} substantially improves \emph{SAA} and generally benefits \emph{Raw}, \emph{KSPCA}, and \emph{RF leaf}. However, \emph{MLKR} and \emph{Neural} perform better with \emph{Obs} than with \emph{E-Rec} in the homoskedastic and conditional-variance DGPs, while they perform similarly in the positive-shock and bimodal settings. Therefore, informative learned representations can already provide strong prescriptions directly from the historical outcomes, without requiring an additional model to estimate the conditional mean. With \emph{O-Rec}, \emph{MLKR} and \emph{Neural} no longer perform worse than with \emph{Obs}, suggesting that the errors in the conditional mean estimates can offset the benefits of recentering.

Under \emph{O-Rec}, the differences between weight constructions nearly disappear in the homoskedastic DGP, and the prescriptions approach the oracle prescription. Here, correcting the conditional mean removes nearly all the systematic differences associated with the features, leaving little for the representation-based weights to exploit. In contrast, differences persist in the conditional-variance DGP and are more pronounced under positive shocks and the bimodal structure, where adjusting the mean does not remove all the differences between the conditional distributions. Representation-based weights can exploit these remaining differences and improve prescriptions even when the scenarios are recentered using the true conditional mean.

\subsection{Effect of the fitting sample size and the feature dimension}
\label{sec:exp4_sample_dimension}

We investigate how prescription performance changes with the fitting sample size and the feature dimension in the linear latent DGP and the nonlinear latent DGP with conditional variance. For the sample-size experiment, we vary the training sample using nested subsets, keeping the 500 validation observations and the SAA and test samples fixed, with \(n=2000\). For the feature-dimension experiment, we increase the number of noise variables while retaining the four demand factors, using 1000 training observations, 500 validation observations, and \(n=1000\). \emph{KSPCA}, \emph{MLKR}, and \emph{Neural} use four-dimensional representations in this comparison.

\begin{figure}[t]
    \centering
    \includegraphics[width=0.8\textwidth]{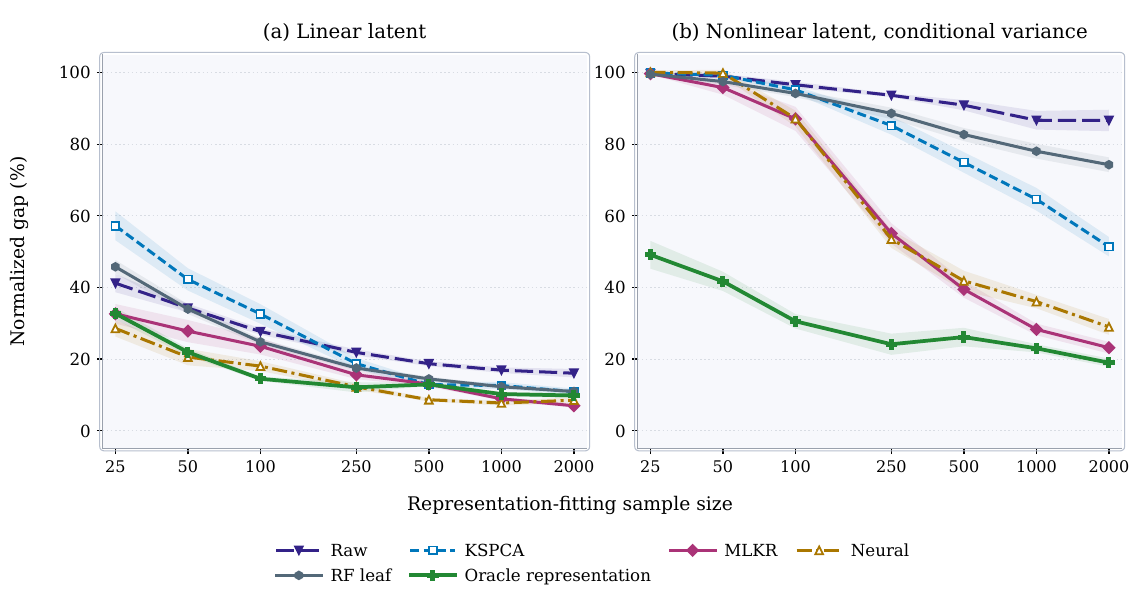}
    \caption{Effect of representation-fitting sample size on the normalized gap.}
    \label{fig:exp4a}
\end{figure}

Figure~\ref{fig:exp4a} shows that the benefit of additional fitting data depends on the DGP. In the linear latent setting, the methods already improve substantially on \emph{SAA} with small training samples. In the nonlinear setting, \emph{MLKR} and \emph{Neural} initially perform close to \emph{SAA} but improve markedly as more training data become available. \emph{KSPCA} improves more gradually, while \emph{Raw} and \emph{RF leaf} retain considerably larger gaps. These results suggest that more flexible learning methods can benefit substantially from additional fitting data, particularly when the relationship between the features and the demand is complex.

\begin{figure}[t]
    \centering
    \includegraphics[width=0.8\textwidth]{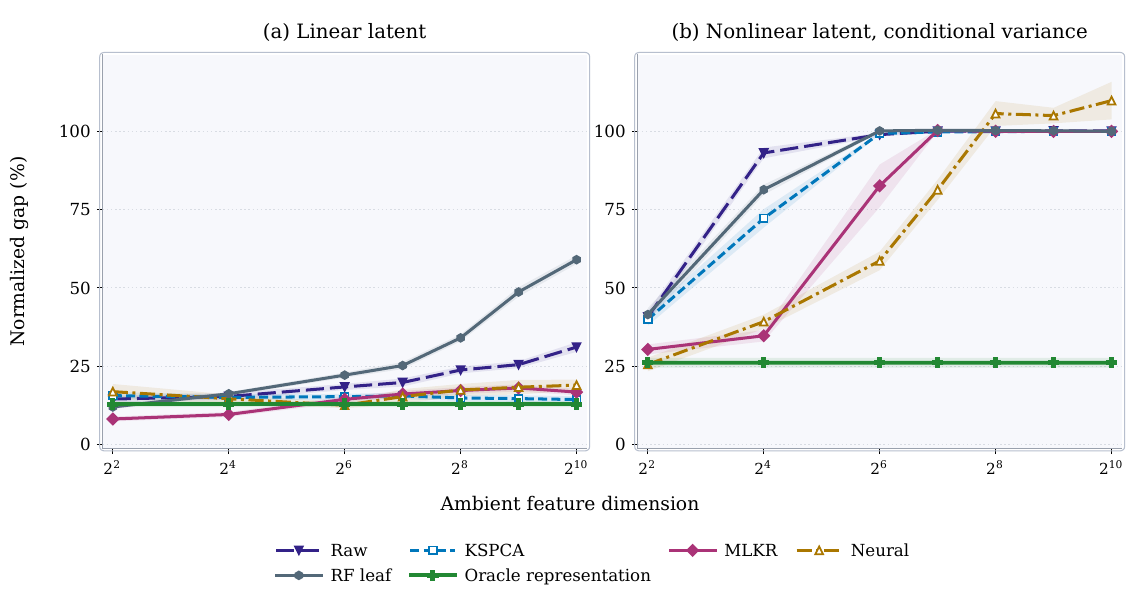}
    \caption{Effect of the feature dimension on the normalized gap.}
    \label{fig:exp4b}
\end{figure}

Figure~\ref{fig:exp4b} shows that adding noise variables affects the methods differently. In the linear latent setting, \emph{KSPCA}, \emph{MLKR}, and \emph{Neural} remain comparatively stable, whereas \emph{Raw} and particularly \emph{RF leaf} deteriorate. This suggests that recovering the relevant demand information becomes more difficult when it is mixed with an increasing number of noise variables, particularly for methods relying on the original features. In the nonlinear setting, all fitted methods eventually approach \emph{SAA}, with \emph{Neural} performing worse than \emph{SAA} at the largest feature dimensions. However, increasing the \emph{Neural} representation dimension to \(d_\phi=32\) brings its performance back to approximately the \emph{SAA} level, without restoring an improvement over \emph{SAA} (see details in Appendix~\ref{app:larger_representations}). These results illustrate the difficulty of learning informative representations from high-dimensional features, even when the underlying demand factors remain unchanged.

\subsection{Effect of additional demand information and representation dimension}
\label{sec:exp_add_info_and_rep_dimension}

Finally, we study how the prescription performance changes as additional demand factors become available and as the dimension of the learned representations increases. For this experiment, we use a separate linear DGP with four demand factors \(U=(U_1,U_2,U_3,U_4)\) of decreasing importance and twelve noise variables, which are mixed through an orthogonal transformation to obtain the observed features. For the oracle benchmarks, we progressively make the first one, two, three, and finally all four demand factors available. For the learned methods, we always use the full 16-dimensional observed feature vector, but fix the dimension of their learned representations to \(d_\phi\in\{1,2,3,4,5,6,8,12,16\}\), rather than tuning it. KSPCA is evaluated only up to \(d_\phi=4\), while Raw uses all 16 observed features independently of \(d_\phi\). We repeat the experiment with SAA sample sizes \(n=100\) and \(n=1000\). Further details of the DGP are provided in Appendix~\ref{app:synthetic_dgps}. For the normalized gap, the oracle prescription used as the reference always has access to all four demand factors \(U\), while the oracle prescription evaluated at each information level uses only the available factors.

\begin{figure}[t]
    \centering
    \includegraphics[width=0.8\textwidth]{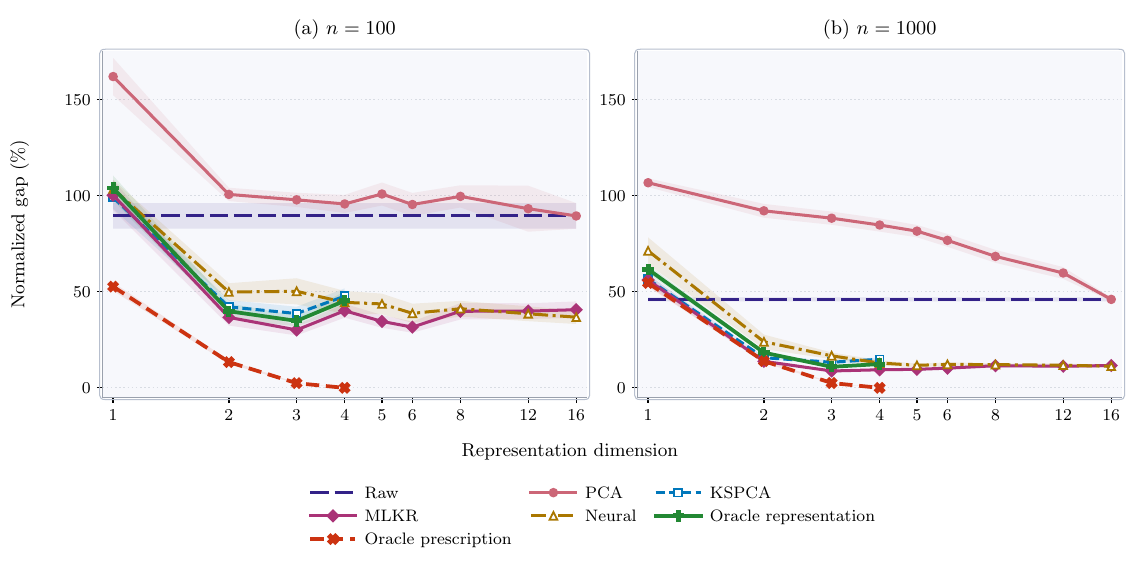}
    \caption{Effect of additional demand information and representation dimension on the normalized gap at \(n=100\) and \(n=1000\). The horizontal axis uses a base-2 logarithmic scale.}
    \label{fig:exp1}
\end{figure}

Figure~\ref{fig:exp1} shows that the oracle benchmarks respond differently as more demand factors become available. The oracle prescription improves with each additional demand factor, although, by construction, \(U_4\) has the smallest contribution to demand. In contrast, the oracle representation improves through three factors but deteriorates when the fourth is added at \(n=100\). With \(n=1000\), its performance moves closer to the oracle prescription, and this deterioration nearly disappears. Thus, additional relevant information may not improve the finite-sample kernel prescription when its contribution to demand is small and the number of historical scenarios is limited. This is because retaining additional factors also increases the dimension in which kernel weights are constructed, reducing the statistical support available for local weighting. This trade-off does not arise for the oracle prescription, which uses the conditional demand distribution directly.

The learned methods also respond differently as representation dimension increases, although they always use the same observed features. At \(n=100\), KSPCA and MLKR improve through three dimensions but deteriorate when the fourth is added, while Neural improves more gradually. However, some methods can recover part of their performance as additional dimensions are retained. For example, MLKR improves again at \(d_\phi=6\), while Neural continues to improve at larger dimensions. Moreover, MLKR and Neural remain substantially better than Raw even when their representation dimension reaches 16. At \(n=1000\), both methods maintain normalized gaps close to 11\%, compared with approximately 46\% for Raw, which uses all 16 observed features throughout the experiment. PCA, in contrast, progressively approaches Raw and coincides with it at \(d_\phi=16\), since retaining all principal components preserves the original Euclidean distances. These results illustrate two potential benefits of learning a feature representation: reducing the dimension in which the scenario weights are constructed and changing how the contexts are compared. The latter can improve the prescription quality even when the dimension of the learned representation is the same as that of the original features.

\subsection{Contextual electricity dispatch}
\label{sec:electricity_dispatch}

We adapt the two-stage electricity dispatch formulation of \citet{zakeri2019pricing} to day-ahead decisions with contextual demand information. Each day, an operator selects hourly generation targets for the following day. After demand is realized, the actual generation can be adjusted upward or downward from these targets, subject to generator capacities and transmission limits, with penalties for unmet demand. The objective includes quadratic production costs and generation adjustment costs. We use the network and production-cost parameters from MATPOWER's \texttt{case30} \citep{zimmerman2011matpower}, scaled to the observed demand level, together with hourly load and temperature data from GEFCom2012 \citep{hong2014gefcom}.

We compare \emph{Raw}, \emph{KSPCA}, \emph{MLKR}, \emph{Neural}, and \emph{RF leaf}, and additionally use representations from Time-Series Mixer (\emph{TSMixer}), a neural forecasting architecture that mixes information across the temporal and feature dimensions \citep{chen2023tsmixer}. Each construction is evaluated using the original historical demand observations and scenarios recentered with predictions from LightGBM or TSMixer. Each historical day receives one scenario weight, which is used in all 24 hourly dispatch problems.

We evaluate the prescriptions over 181 days from  January 1 to June 29, 2008, using a rolling window of the preceding 365 days as the SAA sample. The historical scenarios and their weights are updated daily, while the fitted representations, predictors, and bandwidths remain fixed. We evaluate each prescription by fixing the generation targets and solving the recourse problem for the realized demand. The application also allows us to evaluate the representation-weighted prescriptions in a practical setting that departs from some of the assumptions of the theoretical analysis, as the observations are temporally dependent and the fitting and SAA samples may overlap.

Since the true conditional demand distribution is not known in this application, we cannot compute the oracle prescription. Instead, we use a perfect-information benchmark, obtained by optimizing the generation targets as if the future demand were known in advance. The resulting normalized gap corresponds to the complement of the coefficient of prescriptiveness proposed by \citet{prescriptive_bertsimas}. For the mean daily costs \(C_m\), \(C_{\mathrm{PI}}\), and \(C_{\mathrm{SAA}}\), we report the normalized gap
\begin{equation}
\label{eq:gefcom_normalized_gap}
\mathrm{gap}_m
=
100\,
\frac{C_m-C_{\mathrm{PI}}}
{C_{\mathrm{SAA}}-C_{\mathrm{PI}}}.
\end{equation}
Perfect information corresponds to \(0\%\) and SAA to \(100\%\), with lower values indicating better performance. We also report the mean and the empirical 90th percentile of the daily cost above the perfect-information benchmark. The complete formulation, parameterization, and implementation details are provided in Appendix~\ref{app:gefcom_details}.

\begin{table}[t]
\centering
\caption{Prescription performance in the day-ahead electricity dispatch application.}
\label{tab:gefcom_results}
\footnotesize
\setlength{\tabcolsep}{5pt}
\renewcommand{\arraystretch}{1.10}
\resizebox{\linewidth}{!}{%
\begin{tabular}{@{}lrrrrrrrrr@{}}
\toprule
& \multicolumn{3}{c}{Observations}
& \multicolumn{3}{c}{LightGBM residuals}
& \multicolumn{3}{c}{TSMixer residuals} \\
\cmidrule(lr){2-4}
\cmidrule(lr){5-7}
\cmidrule(lr){8-10}
Method
& Gap (\%) & Mean & \(P_{90}\)
& Gap (\%) & Mean & \(P_{90}\)
& Gap (\%) & Mean & \(P_{90}\) \\
\midrule
SAA
& 100.00 & 60.07 & 104.02
& 17.48 & 10.50 & 18.83
& 15.28 & 9.18 & 15.10 \\
Raw
& 54.37 & 32.66 & 64.47
& 15.70 & 9.43 & 14.24
& 13.99 & 8.40 & 13.26 \\
KSPCA
& 50.99 & 30.63 & 55.22
& 16.17 & 9.71 & 16.85
& 14.97 & 9.00 & 15.12 \\
MLKR
& 45.08 & 27.08 & 45.24
& 14.87 & 8.93 & 14.59
& 13.31 & 7.99 & \textbf{12.47} \\
Neural
& 33.62 & 20.19 & 29.40
& 13.98 & 8.40 & \textbf{12.27}
& \textbf{13.12} & \textbf{7.88} & 12.88 \\
RF leaf
& 30.95 & 18.59 & 29.55
& 14.18 & 8.52 & 13.87
& 13.65 & 8.20 & 13.45 \\
TSMixer
& \textbf{21.55} & \textbf{12.95} & \textbf{17.77}
& \textbf{13.82} & \textbf{8.30} & 13.15
& 13.95 & 8.38 & 14.07 \\
\bottomrule
\end{tabular}%
}
\par\vspace{0.35em}
\begin{minipage}{0.95\linewidth}
\footnotesize
\textit{Notes:} Gap denotes the normalized gap. Mean and \(P_{90}\)
denote the mean and the empirical 90th percentile of daily excess cost
over the same-day perfect-information benchmark, respectively,
in \(10^6\) cost units.
\end{minipage}
\end{table}

Table~\ref{tab:gefcom_results} shows that the choice of representation substantially affects the prescription performance. With the original scenarios, \emph{TSMixer} achieves the lowest mean and 90th-percentile excess costs among the evaluated methods, illustrating the potential of exploiting representations learned by application-specific architectures for contextual prescriptions. However, recentering is particularly effective in this setting: even with equal weights, the recentered scenarios provide a lower mean cost than all other methods using the original scenarios. This contrasts with the synthetic experiments, where informative learned representations sometimes perform better without recentering.

Combining representation-based weights with recentering provides further improvements. \emph{Neural} weighting with TSMixer recentering achieves the lowest mean cost, while \emph{Neural} with LightGBM recentering achieves the lowest 90th-percentile excess cost. These results show that representation-based weighting can complement recentering even when the latter is highly effective. TSMixer further illustrates how the same forecasting architecture can provide both the representations used to construct weights and the predictions used for recentering.

\section{Conclusions}
\label{sec:conclusions}

We developed representation-weighted contextual SAA, which makes the feature representation and the dissimilarity explicit methodological choices in contextual stochastic optimization. By separating these components from the construction of the scenario weights, the framework broadens the range of learning methods available for contextual prescriptions and allows alternative prescriptions to be obtained from the same fitted model, while retaining the estimate-then-optimize structure. The framework also accommodates scenario transformations, allowing representation-based weights to be assigned to adjusted historical outcomes. Our theoretical analysis provides finite-sample bounds on how closely the weighted objective approximates the true conditional expected cost, accounting for the proximity of historical observations, the estimation error in the fitted dissimilarity, sampling variability, and transformation error. For kernel weight constructions, we establish the conditions for obtaining consistent prescriptions.

The synthetic experiments show that the choice of feature representation can substantially affect the prescription quality. \emph{MLKR} and \emph{Neural} achieve particularly strong performance in the nonlinear latent settings, while learned distances and routing-path weights can improve over their conventional constructions. The flexibility of the framework allows decision makers to adapt the choice of the learning method and the representation dimension to the models available to them, the characteristics of their data, and the sizes of the fitting and SAA samples. In particular, the results highlight the benefits of more flexible learning methods when the relationships between the features and the uncertain parameters are complex and sufficient fitting data are available, as well as the potential benefits of retaining fewer dimensions when the number of historical scenarios is limited.

The experiments also highlight the complementary roles of representation-based weighting and scenario transformations. Informative learned representations can provide strong prescriptions without recentering. Even in applications such as the electricity dispatch problem, where recentering is already highly effective, representation-based weighting can further improve the prescription quality. More broadly, decision makers can exploit their existing learning models in different ways without necessarily changing their fitting procedures. A fitted model can provide representations for constructing the scenario weights, predictions for adjusting the historical outcomes, or both. The routing-path results further illustrate how alternative prescriptions can be obtained by exploiting the learned structure of an existing model.

The representation perspective opens several interesting research directions in contextual stochastic optimization. Within the estimate-then-optimize approach, these include exploring alternative representations and dissimilarities from existing learning methods, investigating the prescription quality obtained with other representation-learning methods, and developing learning objectives that capture the aspects of the conditional distribution relevant to the contextual prescriptions. Another direction is to investigate whether small adaptations of the learned representations, for example through fine-tuning, can improve the prescription quality without retraining the entire model. Beyond the estimate-then-optimize approach, an interesting direction is to learn the representations and the dissimilarities using objectives that explicitly account for the downstream optimization problem, connecting the representation perspective with integrated and decision-focused learning approaches.

\section*{Funding}
The first author acknowledges financial support from COLFUTURO (Colombia) through its Loan-Scholarship Program (Crédito Beca).

\section*{Acknowledgments}
Codex (OpenAI) was used for proof auditing, code assistance, and editorial revision. The authors take full responsibility for the manuscript and its results.

\bibliographystyle{plainnat} 
\bibliography{references}

\clearpage
\appendix

\section{Detailed constructions for representation-weighted contextual SAA}
\label{app:framework_examples}

This appendix gives the definitions for the constructions summarized in Table~\ref{tab:constructions}, using the notation of Section~\ref{sec:representation-weighted_SAA}. We leave the dependence on the fitting sample implicit.

\subsection{Basic weighting constructions}
\label{app:basic_weighting}

A kernel prescription in the original feature space uses the identity representation \(\phi(x)=x\), a prescribed dissimilarity such as the Euclidean distance, and the kernel weight construction in \eqref{eq:kernel-weights}.

For nearest-neighbor weights, let \(\mathcal N_{n,k}(x)\) contain the indices of the \(k\) historical feature vectors closest to \(x\) under \(\rho\), with \(1\leq k\leq n\) and the ties resolved by a fixed rule. The weights are
\[
w_i^{\mathrm{knn}}(x)
=
\frac{1}{k}\mathbf 1\{i\in\mathcal N_{n,k}(x)\}.
\]

Membership weights use a group label \(G(x)\) as the representation \(\phi^{\mathrm{mem}}(x)=G(x)\), with the dissimilarity
\[
\rho^{\mathrm{mem}}(x,x')
=
\mathbf 1\{G(x)\neq G(x')\}.
\]
If the group containing \(x\) includes historical observations, the triangular kernel \(K(u)=\max\{1-u,0\}\), with \(h=1\), assigns equal weights to those observations and zero weights to the others. The groups may be specified or estimated, for example through clustering.

\subsection{Tree and ensemble representations}
\label{app:trees_ensembles}

For a fitted tree, let \(L(x)\) denote the terminal leaf reached by \(x\). The representation \(\phi^{\mathrm{leaf}}(x)=L(x)\), combined with the membership weights in Section~\ref{app:basic_weighting}, recovers the conventional tree prescription \citep{prescriptive_bertsimas}.

The routing representation instead records all the nodes visited by \(x\). For the set of tree nodes \(\mathcal V\), define
\[
[\phi^{\mathrm{route}}(x)]_v
=
\mathbf 1\{x\text{ visits node }v\},
\qquad v\in\mathcal V.
\]
The routing-path distance is
\[
\begin{aligned}
\rho^{\mathrm{route}}(x,x')
={}&
\operatorname{depth}(L(x))
+
\operatorname{depth}(L(x'))\\
&-
2\operatorname{depth}\!\left(
\operatorname{LCA}(L(x),L(x'))
\right).
\end{aligned}
\]
Here, the depth of a node is the number of edges from the root, and \(\operatorname{LCA}\) denotes the lowest common ancestor, the last node shared by the two routing paths. This distance is equal to \(\|\phi^{\mathrm{route}}(x)-\phi^{\mathrm{route}}(x')\|_1\) and counts the edges connecting the terminal leaves. Kernel weighting can therefore assign positive weights to observations in different leaves without changing the fitted tree.

A forest combines the weights obtained from its trees. Let \(w_{b,i}(x)\) denote the normalized weight assigned to observation \(i\) by tree \(b\), using either the leaf or routing representation. The forest weights are
\[
w_i^{\mathrm{forest}}(x)
=
\sum_{b=1}^{n_{\mathrm{trees}}}\alpha_b w_{b,i}(x),
\qquad
\alpha_b\geq0,
\qquad
\sum_{b=1}^{n_{\mathrm{trees}}}\alpha_b=1.
\]
The weights from each tree are normalized separately before they are combined. When the same scenarios are used for all the trees, this defines a mixture of their scenario distributions. Equal mixture coefficients \(\alpha_b=1/n_{\mathrm{trees}}\) and uniform weights within each leaf recover the conventional forest construction \citep{prescriptive_bertsimas}.

\subsection{Dimensionality reduction and metric learning}
\label{app:dimensionality_metric}

A linear feature representation has the form
\[
\phi_B(x)=B^\top x,
\]
where \(B\) may be estimated by PCA, supervised PCA, sufficient dimension reduction, or another projection method \citep{dimensionality_ghojogh,spca_barshan,sdr_li}. A nonlinear representation can instead be estimated using kernel PCA or a manifold learning method. The scenario weights are then constructed from the dissimilarities between the represented feature vectors.

Metric learning offers an alternative to prescribing a distance such as the Euclidean or Manhattan distance. It estimates the dissimilarity from data, while the weight construction can remain unchanged. For feature vectors \(x\) and \(x'\), a common construction is the Mahalanobis dissimilarity \citep{metriclearning_bellet},
\[
d_M(x,x')=\sqrt{(x-x')^\top M(x-x')},\qquad M\succeq 0.
\]
Here, \(M\) is a positive semidefinite matrix, and \(M=I\) recovers the Euclidean distance. When \(M=BB^\top\),
\[
d_M(x,x')=\|B^\top(x-x')\|_2.
\]
Thus, the learned dissimilarity can also be computed using the Euclidean distance after a linear transformation. If \(M\) has a low rank, the transformation also reduces the dimension. The same construction applies to a chosen vector representation by replacing \(x\) and \(x'\) with \(\phi(x)\) and \(\phi(x')\). Metric learning for kernel regression (MLKR) learns such a distance by minimizing the leave-one-out kernel-regression prediction error \citep{mlkr_weinberger}.

\subsection{Neural network representations}
\label{app:neural}

For a feed-forward neural network with \(a_0(x)=x\), the activations at layer \(\ell\) satisfy
\[
a_\ell(x)
=
f_\ell\!\left(W_\ell a_{\ell-1}(x)+b_\ell\right),
\]
where \(W_\ell\) and \(b_\ell\) are the fitted parameters and \(f_\ell\) is the activation function. A selected hidden layer \(\ell^\star\) defines the representation
\[
\phi^{\mathrm{neural}}(x)=a_{\ell^\star}(x).
\]
The network may be trained to predict \(\xi\) or through another representation learning objective \citep{representationlearning_bengio}. A prescribed or learned dissimilarity between the selected activations then provides the comparisons used to construct the scenario weights.

\section{Proofs for the statistical guarantees}
\label{app:proofs_statistical_guarantees}
We give the detailed assumptions and proofs for Section~\ref{sec:statistical_guarantees}.

\subsection{Notation used in the proofs}
\label{app:proof_notation}

Fix the observed feature vector \(x\) and write \(w_i=w_{N,n,i}(x)\) for notational simplicity. For the population transformation \(T_x^0\) introduced in Section~\ref{sec:scenario_transformations}, write

\[
 \widetilde\xi_i^0(x)=T_x^0(X_i,\xi_i),\qquad
 \mu_x(z,x')=\mathbb E[c(z,T_x^0(x',\xi))\mid X=x'].
\]
The centered costs used in the sampling argument are
\[
 Y_i(z;x)=c(z,\widetilde\xi_i^0(x))-\mu_x(z,X_i).
\]
All limits involving the fitting and SAA samples are taken jointly as \(N,n\to\infty\).

\subsection{Assumptions}
\begin{assumption}[Independent fitting and SAA samples]
\label{ass:independent_samples}
The fitting sample \(\mathcal D_N^{\mathrm{fit}}\) and the SAA sample \(\mathcal D_n^{\mathrm{saa}}\) are independent i.i.d.\ samples from the joint distribution of \((X,\xi)\). The fitted representation, the dissimilarity, the transformation, and the tuning choices use only the fitting sample. The weights may additionally depend on \(X_1,\ldots,X_n\), but not on the SAA outcomes or the decision \(z\). They are nonnegative and sum to one.
\end{assumption}
\begin{assumption}[Feasible set and cost regularity]
\label{ass:cost_regularity}
The feasible set \(\mathcal Z\subseteq\mathbb R^{d_z}\) is nonempty and compact. There are finite constants \(L_z,L_\xi\ge 0\) such that
\[
 |c(z,\xi)-c(z',\xi')|\le L_z\|z-z'\|+L_\xi\|\xi-\xi'\|
\]
for every \(z,z'\in\mathcal Z\) and all the original and transformed scenario values used in the construction. The conditional expected costs defining \(F(z\mid x)\) and \(\mu_x(z,x')\) are finite.

Both Lipschitz bounds hold, for example, in two-stage linear programs with fixed recourse, right-hand-side uncertainty, and finite-valued recourse functions. They also hold for more general convex recourse problems with uniformly bounded subgradients and a Lipschitz first-stage cost.

\end{assumption}
\begin{assumption}[Conditional tails]
\label{ass:subgaussian}
Conditional on the fitting sample and \(X_1,\ldots,X_n\), the centered costs \(Y_i(z;x)\) are sub-Gaussian with a variance proxy at most \(\bar\tau_x^2<\infty\), uniformly over \(i=1,\ldots,n\) and \(z\in\mathcal Z\). Specifically, for every \(\lambda\in\mathbb R\),
\[
 \mathbb E\!\left[e^{\lambda Y_i(z;x)}\mid\mathcal D_N^{\mathrm{fit}},X_1,\ldots,X_n\right]
 \le e^{\lambda^2\bar\tau_x^2/2}.
\]
The constant \(\bar\tau_x\) does not depend on \(N,n\). A sub-Gaussian random variable has tails that decay at a Gaussian rate or faster. In particular, uniformly bounded population scenario costs satisfy this condition.
\end{assumption}
\begin{assumption}[Relevance of the population dissimilarity]
\label{ass:relevance_dissimilarity}
For the chosen population transformation, there is a finite constant \(L_\rho\ge 0\) such that
\[
 \Gamma_x^T(x')\le L_\rho\rho_0(x,x'),\qquad x'\in\mathcal X.
\]
\end{assumption}
For the original scenarios, Assumption~\ref{ass:relevance_dissimilarity} corresponds to condition~\ref{eq:population_relevance_condition}. For the transformed scenarios, it uses \(\Gamma_x^T\) in place of \(\Gamma_x\). The population dissimilarity therefore controls the conditional cost differences that remain after the population transformation.

\subsection{Proof of the bound using the fitted representation}
\begin{proof}[Proof of Theorem~\ref{cor:fitted_representation_bound}]
Condition on the fitting sample and the SAA feature vectors. Under Assumption~\ref{ass:independent_samples}, the weights are fixed and the centered costs \(Y_i(z;x)\) are independent and have mean zero for each fixed decision \(z\).

\emph{Separating the errors.} For each feasible decision, adding and subtracting the weighted population scenario costs and their conditional expectations gives
\begin{equation}\label{eq:lfrB-decomp}
\begin{aligned}
 \widehat F_{N,n}(z\mid x)-F(z\mid x)
 &=\sum_{i=1}^n w_i\big[c(z,\widetilde\xi_{N,i}(x))-c(z,\widetilde\xi_i^0(x))\big]\\
 &\quad+\sum_{i=1}^n w_i\big[\mu_x(z,X_i)-F(z\mid x)\big]\\
 &\quad+\sum_{i=1}^n w_iY_i(z;x),
\end{aligned}
\end{equation}
where we use $\sum_i w_i=1$. The three terms distinguish the error in estimating the transformation, the differences in conditional expected costs that remain even with the population transformation, and the sampling variability. We bound each term uniformly over the feasible decisions.

\emph{The transformation error.} The largest absolute value of the first term is the transformation error from Section~\ref{sec:scenario_transformations}:
\begin{equation}\label{eq:lfrB-transform}
 E^T_{N,n}(x)=\sup_{z\in\mathcal Z}\left|\sum_{i=1}^n w_i\big[c(z,\widetilde\xi_{N,i}(x))-c(z,\widetilde\xi_i^0(x))\big]\right|.
\end{equation}
By the Lipschitz bound in Assumption~\ref{ass:cost_regularity},
\begin{equation}\label{eq:lfrB-transform-lip}
 E^T_{N,n}(x)\le L_\xi\sum_{i=1}^n w_i\|\widetilde\xi_{N,i}(x)-\widetilde\xi_i^0(x)\|.
\end{equation}
Without a transformation, this error is zero.

\emph{The weighting bias.} The absolute difference at each \(X_i\) is at most \(\Gamma_x^T(X_i)\), whatever decision we evaluate. Assumption~\ref{ass:relevance_dissimilarity} lets us apply the weighted-radius argument of Section~\ref{sec:role-dissimilarity} to these remaining discrepancies, with \(\Gamma_x^T\) in place of \(\Gamma_x\). Thus the second term in~\eqref{eq:lfrB-decomp} satisfies
\begin{equation}\label{eq:lfrB-cost}
\begin{aligned}
 \sup_{z\in\mathcal Z}\left|\sum_{i=1}^n w_i\big[\mu_x(z,X_i)-F(z\mid x)\big]\right|
 &\le \sum_{i=1}^n w_i\Gamma_x^T(X_i)\\
 &\le L_\rho\sum_{i=1}^n w_i\rho_0(x,X_i)\\
 &\le L_\rho\big(\widehat R_{N,n}(x)+\Delta_{N,n}(x)\big).
\end{aligned}
\end{equation}

\emph{The sampling variability.} For each fixed decision, conditional independence and Assumption~\ref{ass:subgaussian} make \(\sum_i w_iY_i(z;x)\) sub-Gaussian with a variance proxy at most \(\bar\tau_x^2\sum_i w_i^2\). The squared weights arise because multiplying a centered cost by \(w_i\) multiplies its variance proxy by \(w_i^2\).

For \(\eta>0\), take a smallest \(\eta\)-cover \(z_1,\ldots,z_J\) of \(\mathcal Z\), where \(J=\mathcal N_{\mathcal Z}(\eta)\). Every feasible decision lies within distance \(\eta\) of a cover point, and compactness ensures that \(J\) is finite. Applying Hoeffding's inequality for independent sub-Gaussian variables at each cover point with failure probability \(\delta/J\), then the union bound, gives, with conditional probability at least \(1-\delta\),
\[
 \max_{1\le j\le J}\left|\sum_{i=1}^n w_iY_i(z_j;x)\right|
 \le\bar\tau_x\sqrt{2\log(2J/\delta)\sum_{i=1}^n w_i^2}.
\]
Both the population scenario cost and its conditional expectation are \(L_z\)-Lipschitz in the decision. Their difference, and hence its weighted average, is therefore \(2L_z\)-Lipschitz. Moving from a cover point to a decision at distance at most \(\eta\) adds at most \(2L_z\eta\), so
\[
 \sup_{z\in\mathcal Z}\left|\sum_{i=1}^n w_iY_i(z;x)\right|
 \le\max_{1\le j\le J}\left|\sum_{i=1}^n w_iY_i(z_j;x)\right|+2L_z\eta.
\]
This is the standard SAA argument using finite covers \citep[Chapter~5]{saa_shapiro}, with \(\sum_i w_i^2\) replacing \(1/n\). Taking the supremum in~\eqref{eq:lfrB-decomp} and combining the three bounds by the triangle inequality gives~\ref{eq:fitted_representation_bound}.
\end{proof}

\subsection{Kernel weights and consistency}
\label{app:kernel_consistency}
We use the following condition on the probability of finding SAA observations near \(x\) to establish the kernel bounds needed for Corollary~\ref{cor:kernel_consistency}.

\begin{assumption}[Probability of nearby observations]
\label{ass:local_mass}
There are \(q>0\), \(s_0>0\), and \(c_->0\) such that
\[
 \mathbb P\{\rho_0(x,X)\le s\}\ge c_-s^q,\qquad 0<s\le s_0.
\]
The exponent \(q\) and the constants may depend on \(x\), but not on the sample sizes.
\end{assumption}

We retain the sampling conditions of Assumption~\ref{ass:independent_samples} and select the bandwidth using the fitting sample. If the kernel denominator is zero, we assign equal weights \(w_i=1/n\).

We use the following proposition to show that the fitted weighted radius tends to zero even as the effective sample size grows. The idea is to consider observations within a small fraction of the population bandwidth, leaving room for the estimation error in their fitted distances. Each such observation then contributes at least a fixed positive amount to the kernel denominator. By Assumption~\ref{ass:local_mass}, their conditional expected number is at least a constant times \(nh_N^q\), so their combined contribution grows as \(nh_N^q\xrightarrow{p}\infty\).

\begin{proposition}[Kernel weights under a fitted dissimilarity]\label{prop:lfrB-kernel}
Under these sampling and bandwidth conditions, suppose Assumption~\ref{ass:local_mass} holds and
\[
 h_N\xrightarrow{p}0,\qquad nh_N^q\xrightarrow{p}\infty,\qquad \Delta_{N,n}(x)/h_N\xrightarrow{p}0.
\]
The nonnegative kernel \(K\) is either bounded, continuous at zero, positive at zero, and zero for \(u>1\), or is the Gaussian kernel \(K(u)=\exp(-u^2/2)\). There is a constant \(c_{k_0}>0\), independent of the sample sizes, such that
\begin{equation}\label{eq:lfrB-kernel}
 \mathbb P\{n_{\mathrm{eff},N,n}(x)\ge c_{k_0}nh_N^q\}\longrightarrow 1,
 \qquad \widehat R_{N,n}(x)\xrightarrow{p}0.
\end{equation}
For kernels with compact support, the probability of a zero kernel denominator tends to zero.
\end{proposition}

\begin{proof}
Write \(h=h_N\), \(d_i=\rho_0(x,X_i)\), \(r_i=\widehat\rho_N(x,X_i)\), and \(\Delta=\Delta_{N,n}(x)\). Let \(K_i=K(r_i/h)\) and \(S_n=\sum_i K_i\). Multiplying a bounded kernel by a positive constant does not change the normalized weights, so we may take \(0\le K_i\le 1\). Conditional on the fitting sample, \(h\) is fixed and the SAA feature vectors remain independent.

\emph{The effective sample size.} Since \(K_i^2\le K_i\), whenever \(S_n>0\) we have
\[
 n_{\mathrm{eff},N,n}(x)=\frac{1}{\sum_i w_i^2}
 =\frac{S_n^2}{\sum_i K_i^2}\ge\frac{S_n^2}{\sum_i K_i}=S_n.
\]
Both kernel classes are continuous at zero and satisfy \(K(0)>0\). Choose \(a\in(0,1)\) and \(k_0>0\) such that \(K(u)\ge k_0\) for \(0\le u\le a\). Since \(\Delta/h\xrightarrow{p}0\), we have \(\Delta\le ah/2\) with probability tending to one. Whenever this inequality holds,
\[
 d_i\le ah/2 \quad\Longrightarrow\quad r_i\le d_i+\Delta\le ah \quad\Longrightarrow\quad K_i\ge k_0.
\]
Thus \(S_n\ge k_0\sum_{i=1}^n\mathbf{1}\{d_i\le ah/2\}\) with probability tending to one.

Conditional on the fitting sample, the number of observations with \(d_i\le ah/2\) is binomial. Whenever \(ah/2\le s_0\), its conditional mean is at least \(nc_-(ah/2)^q\) by Assumption~\ref{ass:local_mass}. Since \(h\xrightarrow{p}0\) and \(nh^q\xrightarrow{p}\infty\), binomial concentration gives, for some constant \(c_{k_0}>0\) independent of \(N,n\), with probability tending to one,
\begin{equation}
\label{eq:lfrB-denom}
n_{\mathrm{eff},N,n}(x)\ge S_n\ge c_{k_0}nh^q.
\end{equation}
In particular, \(\sum_i w_i^2=1/n_{\mathrm{eff},N,n}(x)\xrightarrow{p}0\). Since \(\mathbb P\{S_n=0\}\to0\), the choice of weights \(w_i=1/n\) when the denominator is zero does not affect these limits.

\emph{The fitted weighted radius.} When \(S_n>0\), the weighted average fitted distance is
\[
 \widehat R_{N,n}(x)=\frac{\sum_{i=1}^n r_iK_i}{S_n}.
\]
For a compactly supported kernel, only observations with \(r_i\le h\) contribute to the numerator. Their weighted average is therefore at most \(h\). Since \(h\xrightarrow{p}0\), it follows that \(\widehat R_{N,n}(x)\xrightarrow{p}0\).

The Gaussian kernel has no such cutoff: even observations with \(r_i>h\) receive positive weight. However, their contribution to the numerator, \(r_i\exp(-r_i^2/(2h^2))\), reaches its maximum at \(r_i=h\) and decreases thereafter. To show that the weighted radius vanishes, fix \(\varepsilon>0\). Since \(h\xrightarrow{p}0\), we have \(h<\varepsilon/2\) with probability tending to one. Observations with \(r_i\le\varepsilon/2\) contribute at most \(\varepsilon/2\) to the weighted radius because their weights sum to at most one. For the remaining observations, the decreasing Gaussian contribution is bounded by its value at \(\varepsilon/2\). There are at most \(n\) such observations. Applying the denominator bound~\eqref{eq:lfrB-denom} gives, with probability tending to one,
\[
\begin{aligned}
\widehat R_{N,n}(x)
&\le \frac{\varepsilon}{2}
+\frac{1}{S_n}\sum_{i:r_i>\varepsilon/2}
r_i\exp\!\left(-\frac{r_i^2}{2h^2}\right)\\
&\le \frac{\varepsilon}{2}
+\frac{n}{S_n}\frac{\varepsilon}{2}
\exp\!\left(-\frac{\varepsilon^2}{8h^2}\right)\\
&\le \frac{\varepsilon}{2}
+\frac{\varepsilon}{2c_{k_0}h^q}
\exp\!\left(-\frac{\varepsilon^2}{8h^2}\right).
\end{aligned}
\]
Since \(h\xrightarrow{p}0\), the exponential decay dominates \(h^{-q}\), so the second term converges to zero in probability. Hence \(\mathbb P\{\widehat R_{N,n}(x)>\varepsilon\}\to0\) for every \(\varepsilon>0\).
\end{proof}

\begin{proof}[Proof of Corollary~\ref{cor:kernel_consistency}]
Under the stated conditions, Proposition~\ref{prop:lfrB-kernel} gives
\[
n_{\mathrm{eff},N,n}(x)\xrightarrow{p}\infty,\qquad \widehat R_{N,n}(x)\xrightarrow{p}0.
\]
Since \(\Delta_{N,n}(x)/h_N\xrightarrow{p}0\) and \(h_N\xrightarrow{p}0\), we also have \(\Delta_{N,n}(x)\xrightarrow{p}0\). The transformation error \(E^T_{N,n}(x)\) converges to zero in probability by assumption.

Fix \(\varepsilon>0\) and \(0<\delta<1\). Choose a cover radius \(\eta>0\) such that \(2L_z\eta\le\varepsilon/2\). By compactness of \(\mathcal Z\), the cover size \(J=\mathcal N_{\mathcal Z}(\eta)\) is finite and does not depend on the sample sizes. Applying Theorem~\ref{cor:fitted_representation_bound} with failure probability \(\delta/2\), the sampling term is
\[
\bar\tau_x\sqrt{\frac{2\log(4J/\delta)}{n_{\mathrm{eff},N,n}(x)}}\xrightarrow{p}0,
\]
because the numerator is fixed and the effective sample size diverges.

Let \(B_{N,n}\) be the sum of all terms in the bound except the covering error. Each of these terms converges to zero in probability, so \(B_{N,n}\xrightarrow{p}0\). Since the covering error is at most \(\varepsilon/2\), the uniform objective error can exceed \(\varepsilon\) only if the bound fails or \(B_{N,n}>\varepsilon/2\). Hence, for all sufficiently large \(N,n\),
\[
\mathbb P\!\left\{\sup_{z\in\mathcal Z}\bigl|\widehat F_{N,n}(z\mid x)-F(z\mid x)\bigr|>\varepsilon\right\}
\le \frac{\delta}{2}+\mathbb P\{B_{N,n}>\varepsilon/2\}\le\delta.
\]
As \(\varepsilon\) and \(\delta\) were arbitrary, this proves
\begin{equation}\label{eq:lfrB-uniform}
\sup_{z\in\mathcal Z}\bigl|\widehat F_{N,n}(z\mid x)-F(z\mid x)\bigr|\xrightarrow{p}0.
\end{equation}

Prescription consistency then follows directly from \eqref{eq:lfrB-regret}.

\end{proof}
\section{Verification for selected constructions}
\label{app:method_verification}

\subsection{Spectral representations}
\label{app:spectral_verification}

We verify the dissimilarity estimation condition \(\Delta_{N,n}(x)/h_N\xrightarrow{p}0\) used in Corollary~\ref{cor:kernel_consistency} for PCA, supervised PCA (SPCA), kernel PCA (KPCA), and kernel supervised PCA (KSPCA). These methods project the features onto a subspace determined by eigenvectors. The common argument is that accurately estimating this subspace also gives accurate distances between the projected feature vectors. We follow the constructions in \citet[Chapter~5]{dimensionality_ghojogh}, using a linear response kernel for the supervised methods. The representation kernel may be linear or nonlinear. All the estimates use the independent fitting sample \(\mathcal D_N^{\mathrm{fit}}\). The retained dimension \(d_\phi\), the response dimension \(d_\xi\), and the representation kernel are fixed. The original feature dimension is also fixed for the linear methods.

Let \(\Phi:\mathcal X\to\mathcal H\) be the feature map associated with the positive semidefinite representation kernel \(\kappa\), distinct from the weighting kernel \(K\), and assume
\[\sup_{u\in\mathcal X}\|\Phi(u)\|_{\mathcal H}^2=\sup_{u\in\mathcal X}\kappa(u,u)\leq\kappa_{\max}<\infty.\]
For the linear kernel, \(\Phi(u)=u\), so this is a boundedness assumption on the original features.

For each method, let \(A_0\) denote the population operator and \(\widehat A_N\) its estimate. These are matrices for the linear methods, whereas for the kernel methods they act on a feature space \(\mathcal H\) that may be infinite-dimensional. Let \(P_0\) and \(\widehat P_N\) be the orthogonal projectors onto their respective leading \(d_\phi\)-dimensional subspaces. The corresponding representations are \(\phi_0(u)=P_0\Phi(u)\) and \(\widehat\phi_N(u)=\widehat P_N\Phi(u)\). Distances between projected features depend on the retained subspace, not on its choice of orthonormal basis. Centering also leaves pairwise differences unchanged. By the reverse triangle inequality and boundedness of \(\Phi\),
\begin{equation}\Delta_{N,n}(x)\le 2\sqrt{\kappa_{\max}}\,\|\widehat P_N-P_0\|_{\mathrm{op}}, \label{eq:projector_dissimilarity_bridge}\end{equation}
where the operator norm \(\|A\|_{\mathrm{op}}=\sup_{\|v\|=1}\|Av\|\) measures how much \(A\) can stretch a unit vector.

To bound the effect of estimating \(A_0\) on the retained subspace, we follow the argument given for matrices by \citet[Section~2.3.4, proof of Corollary~2.8]{spectralmethods_chen}. Set \(r=d_\phi\) and let \(\lambda_j(A)\) denote the \(j\)-th largest eigenvalue of \(A\), counting multiplicities. Assume a positive population spectral gap \(g=\lambda_r(A_0)-\lambda_{r+1}(A_0)>0\) between the last retained and first discarded eigenvalues. Whenever \(\|\widehat A_N-A_0\|_{\mathrm{op}}<g/2\), we have
\[
\begin{aligned}
\|\widehat P_N-P_0\|_{\mathrm{op}}
&\le \frac{\|\widehat A_N-A_0\|_{\mathrm{op}}}{\lambda_r(A_0)-\lambda_{r+1}(\widehat A_N)}\\
&\le \frac{\|\widehat A_N-A_0\|_{\mathrm{op}}}{g-\|\widehat A_N-A_0\|_{\mathrm{op}}}
\le \frac{2\|\widehat A_N-A_0\|_{\mathrm{op}}}{g}.
\end{aligned}
\]
The first inequality is the \(\sin\Theta\) bound of \citet[eq.~(6.1)]{rotationeigen_Kahan}: it bounds the rotation of the retained subspace using the separation between the retained population eigenvalues and the discarded fitted eigenvalues. The second follows from Weyl's inequality: the first discarded fitted eigenvalue can exceed its population counterpart by at most \(\|\widehat A_N-A_0\|_{\mathrm{op}}\). The last inequality uses the fact that at least half the population spectral gap remains. For matrices, \citet[Exercise~4.13]{vershynin2026high} states the final bound directly.

For PCA and KPCA, \(\widehat A_N\) estimates \(A_0=\operatorname{Cov}(\Phi(X))\) by averaging products of centered features over the \(N\) independent fitting observations. Under the boundedness assumption on \(\Phi\), the covariance estimation error \(\|\widehat A_N-A_0\|_{\mathrm{op}}\) has the usual \(1/\sqrt N\) sampling scale. This follows from the Hilbert--Schmidt norm bound of \citet[Lemma~1 and the discussion of centering in Section~3]{kpca_zwald}, which also bounds the operator norm.

For SPCA and KSPCA, the construction of \citet{spca_barshan} with the linear response kernel \(\ell(\xi,\xi')=\xi^\top\xi'\) uses \(A_0=CC^*\), where \(C=\operatorname{Cov}(\Phi(X),\xi)\) and \({}^*\) denotes the adjoint. Assume also that \(\xi\) is bounded. Since \(C\) is a block of the covariance of \((\Phi(X),\xi)\), the same sampling bound applies to its estimate \(\widehat C_N\). It also applies to \(\widehat A_N=\widehat C_N\widehat C_N^*\), since
\[
\|\widehat A_N-A_0\|_{\mathrm{op}}\le \bigl(\|\widehat C_N\|_{\mathrm{op}}+\|C\|_{\mathrm{op}}\bigr)\|\widehat C_N-C\|_{\mathrm{op}}
\]
and the first factor is bounded. For both supervised methods, \(\operatorname{rank}(A_0)\le d_\xi\), so the eigengap condition requires \(d_\phi\le d_\xi\).

With a fixed positive eigengap, the projector bound and \eqref{eq:projector_dissimilarity_bridge} give the same \(1/\sqrt N\) sampling bound for the error in the fitted distances. Choosing \(h_N\) to decrease more slowly than \(1/\sqrt N\) therefore gives \(\Delta_{N,n}(x)/h_N\xrightarrow{p}0\).

\subsection{Estimated Mahalanobis dissimilarities}
\label{app:mahalanobis_verification}

We bound the dissimilarity error \(\Delta_{N,n}(x)\) in terms of the error in the fitted Mahalanobis matrix. Let \(M_0,\widehat M_N\in\mathbb R^{d_x\times d_x}\) be positive semidefinite matrices, where \(M_0\) is fixed and \(\widehat M_N\) is estimated from the fitting sample \(\mathcal D_N^{\mathrm{fit}}\). The population and fitted dissimilarities are
\[
\rho_0(x,x')=\sqrt{(x-x')^\top M_0(x-x')},\qquad \widehat\rho_N(x,x')=\sqrt{(x-x')^\top\widehat M_N(x-x')}.
\]

Write \(e=\|\widehat M_N-M_0\|_{\mathrm{op}}\) and assume \(C_X=\sup_{u\in\mathcal X}\|u\|_2<\infty\). The matrix error directly controls the difference between the squared dissimilarities:
\[
\bigl|\widehat\rho_N(x,x')^2-\rho_0(x,x')^2\bigr|=\bigl|(x-x')^\top(\widehat M_N-M_0)(x-x')\bigr|\le e\,\|x-x'\|_2^2.
\]
Since \(|\sqrt a-\sqrt b|\le\sqrt{|a-b|}\) for \(a,b\ge0\) and \(\|x-x'\|_2\le 2C_X\), taking the maximum over \(x'=X_1,\ldots,X_n\) gives
\begin{equation}
\label{eq:mahalanobis-distance-error}
\Delta_{N,n}(x)\le 2C_X\sqrt e.
\end{equation}
Thus \(e/h_N^2\xrightarrow{p}0\) is sufficient for the dissimilarity estimation condition \(\Delta_{N,n}(x)/h_N\xrightarrow{p}0\) in Corollary~\ref{cor:kernel_consistency}.

Our experiments use multi-output MLKR, which fits a shared Mahalanobis matrix by minimizing the sum of leave-one-out squared prediction errors across outcome components. The bound above applies regardless of how the matrix is fitted. However, different matrices can yield equally accurate conditional mean predictions while assigning different distances to the same contexts. Establishing when multi-output MLKR estimates \(M_0\) accurately enough to satisfy the bandwidth condition remains a direction for future research.

\subsection{Tree and forest weights}
\label{app:tree_consistency}
Consistency results for tree prescriptions usually restrict how trees split and grow so that the target leaf covers a smaller region of the feature space while its number of observations grows (\citet[Section~3.2]{prescriptivetrees_bertsimasmccord}, \citet[Section~5]{stochasticforest_kallus}, and \citet[Section~2]{honesttrees_bladt}).

\begin{figure}[!htbp]
\centering
\begin{subcaptiongroup}%
\phantomcaption\label{fig:tree-ancestor-region}%
\includegraphics[height=3.3cm]{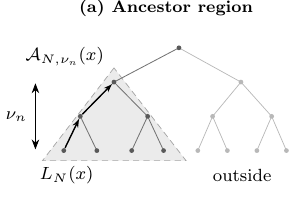}\hspace{0.55cm}%
\phantomcaption\label{fig:tree-growing-regions}%
\includegraphics[height=3.3cm]{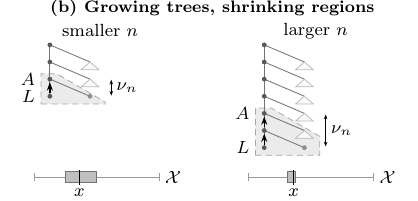}%
\end{subcaptiongroup}
\captionsetup{font=small,skip=3pt}
\caption{(\subref*{fig:tree-ancestor-region}) The shaded subtree contains all the leaves below the selected ancestor. (\subref*{fig:tree-growing-regions}) The target leaf \(L\) is at depth 3 in the first tree and 6 in the second. Moving back 1 and 2 levels, respectively, places the ancestor \(A\) at depths 2 and 4. The shaded intervals show the corresponding regions containing \(x\).}\label{fig:tree-localization}
\end{figure}

In our setting, conventional leaf weights use only observations in the same leaf as \(x\), whereas Gaussian routing weights also use observations from other leaves. Our proof separates the observations inside the region of a selected ancestor from those outside it (see Figure~\ref{fig:tree-localization}(\subref{fig:tree-ancestor-region})). Inside, the observations together receive nearly all the weight, so their contexts must become close to \(x\). Outside, the contexts can be far from \(x\), but their total weight vanishes.

Fix \(x\in\mathcal X\). Retain Assumptions~\ref{ass:independent_samples}--\ref{ass:relevance_dissimilarity} for the chosen scenario construction and suppose \(E^T_{N,n}(x)\xrightarrow{p}0\). The tree and its bandwidth \(0<h_N<\infty\) are determined independently of the SAA observations, using the fitting sample, the prescribed size \(n\), and any independent tree randomization. Unlike the \hyperref[app:spectral_verification]{spectral results}, this proof does not require the fitted routing distance to estimate \(\rho_0\). Since the rule at each internal node compares one of the original features with a threshold, we additionally assume for the construction below that closeness in the original feature space implies closeness under the population dissimilarity, that is, \(\rho_0(x,x')\to0\) as \(x'\to x\).

\paragraph{Routing-path weights.}
With a slight abuse of notation, let \(L_N(x)\) denote the set of contexts assigned to the same terminal leaf as \(x\). Let \(\mathcal I_{N,n}(x)=\{i\in\{1,\ldots,n\}:X_i\in L_N(x)\}\) index the SAA observations in that leaf, with \(|\mathcal I_{N,n}(x)|\) observations in total. The routing distance \(\widehat\rho_N^{\mathrm{route}}(x,x')\) counts the edges on the path between the leaves of \(x\) and \(x'\), and it is zero when both contexts reach the same leaf.

Starting at the target leaf, move \(\nu\) edges toward the root, stopping at the root if necessary. Let \(\mathcal A_{N,\nu}(x)\) be the union of the leaf regions below the selected ancestor, as shown in Figure~\ref{fig:tree-localization}(\subref{fig:tree-ancestor-region}). To reach a leaf outside this region, the path from the target leaf must pass the selected ancestor and go at least one edge farther toward the root, so its length exceeds \(\nu\).

We choose \(\nu_n\) so that each outside observation contributes at most \(1/n\) before normalization. With \(n\) observations in total, their combined contribution is then at most one. For the Gaussian kernel, take
\[
\nu_n=\big\lceil h_N\sqrt{2\log n}\big\rceil,
\qquad \exp\!\left(-\frac{\nu_n^2}{2h_N^2}\right)\le \frac1n.
\]

\begin{proposition}[Consistency of routing weights]
\label{prop:routing_consistency}
In addition to the preceding assumptions, assume that the population dissimilarities from \(x\) are bounded over the entire context space: \(\bar\rho_0(x):=\sup_{x'\in\mathcal X}\rho_0(x,x')<\infty\). For an independent context \(X\), let \(p_N(x)=\mathbb P\{X\in L_N(x)\mid\text{fitted tree}\}\). Suppose that, as \(N,n\to\infty\),
\begin{gather}
 n p_N(x)\xrightarrow{p}\infty,\label{eq:tree-leaf-count}\\
 \Lambda_{N,n}(x):=\sup_{x'\in\mathcal A_{N,\nu_n}(x)}\rho_0(x,x')\xrightarrow{p}0.\label{eq:tree-ancestor-shrinkage}
\end{gather}
Here \(\Lambda_{N,n}(x)\) measures the greatest population dissimilarity from \(x\) over the region defined by the selected ancestor. Then the Gaussian routing weights satisfy
\begin{equation}\label{eq:tree-objective-consistency}
 \sup_{z\in\mathcal Z}|\widehat F_{N,n}(z\mid x)-F(z\mid x)|\xrightarrow{p}0.
\end{equation}
\end{proposition}

\begin{proof}
\textit{Effective sample size.} Conditional on the fitted tree, \(|\mathcal I_{N,n}(x)|\) is binomial with mean \(np_N(x)\), so \eqref{eq:tree-leaf-count} gives \(|\mathcal I_{N,n}(x)|\xrightarrow{p}\infty\). Now consider the weights assigned to observations in the target leaf. Each contributes \(\exp(0)=1\) to the sum of the Gaussian kernel values. With
\[
 K_i:=\exp\!\left(-\frac{[\widehat\rho_N^{\mathrm{route}}(x,X_i)]^2}{2h_N^2}\right),
 \qquad S_n:=\sum_{i=1}^n K_i,\qquad w_{N,n,i}(x):=K_i/S_n,
\]
we therefore have \(S_n\ge |\mathcal I_{N,n}(x)|\xrightarrow{p}\infty\). As in the \hyperref[prop:lfrB-kernel]{general kernel proof}, \(K_i^2\le K_i\), so
\[
 n_{\mathrm{eff},N,n}(x)=\frac{S_n^2}{\sum_iK_i^2}\ge S_n\xrightarrow{p}\infty.
\]

\textit{Weighted population dissimilarity.} Now consider the weights outside the ancestor region \(\mathcal A=\mathcal A_{N,\nu_n}(x)\). Here we use the additional bound \(\bar\rho_0(x)<\infty\). Since the outside kernel values sum to at most one, the full outside contribution satisfies
\[
\begin{aligned}
 \sum_{i:X_i\notin\mathcal A}w_{N,n,i}(x)\rho_0(x,X_i)
 &=\frac{1}{S_n}\sum_{i:X_i\notin\mathcal A}K_i\rho_0(x,X_i)\\
 &\le \frac{\bar\rho_0(x)}{S_n}\sum_{i:X_i\notin\mathcal A}K_i\\
 &\le \frac{\bar\rho_0(x)}{S_n}\xrightarrow{p}0.
\end{aligned}
\]
Inside the ancestor region, each dissimilarity is at most \(\Lambda_{N,n}(x)\) and the weights sum to at most one. By \eqref{eq:tree-ancestor-shrinkage},
\[
\begin{aligned}
 \sum_{i:X_i\in\mathcal A} w_{N,n,i}(x)\rho_0(x,X_i)
 &\le \Lambda_{N,n}(x)\sum_{i:X_i\in\mathcal A}w_{N,n,i}(x)\\
 &\le \Lambda_{N,n}(x)\xrightarrow{p}0.
\end{aligned}
\]
Thus both \(\sum_i w_{N,n,i}(x)^2=1/n_{\mathrm{eff},N,n}(x)\) and \(\sum_i w_{N,n,i}(x)\rho_0(x,X_i)\) converge to zero in probability. Together with \(E^T_{N,n}(x)\xrightarrow{p}0\), these limits yield \eqref{eq:tree-objective-consistency} by the same uniform concentration argument used in the proof of Theorem~\ref{cor:fitted_representation_bound}.
\end{proof}

\paragraph{A sufficient tree construction.}
We give one construction that satisfies \eqref{eq:tree-leaf-count} and \eqref{eq:tree-ancestor-shrinkage}. It must retain enough probability in the target leaf for its expected SAA count to grow, while bringing all contexts below the selected ancestor close to \(x\). The tree splits on the observed features, so the split rules control how much each feature can vary below the selected ancestor. As these ranges narrow around the feature values of \(x\), the continuity condition above ensures that every context in the region becomes close to \(x\) under \(\rho_0\). Assumption~\ref{ass:relevance_dissimilarity} then bounds the corresponding differences in conditional expected costs. Assume also that \(\mathcal X=[0,1]^{d_x}\) and that the context density is bounded below by \(m>0\) on this cube.

Cycle through the features in a fixed order. At each node, use the regression split rule of classification and regression trees (CART) to minimize the sum of within-child squared errors on the fitting data, subject to each child retaining at least a fraction \(\theta\in(0,1/2)\) of the selected feature's range (\citet[Definition~3 and Algorithm~1]{regulartreesearch_wang}). Use the midpoint when no admissible candidate separates the fitting observations. Grow every branch to depth \(d_n=\lfloor\log n/(2\log(1/\theta))\rfloor\). This depth keeps the target leaf's volume at least \(\theta^{d_n}\ge n^{-1/2}\), so its expected SAA count is at least \(m\sqrt n\), giving \eqref{eq:tree-leaf-count}.

For this construction, take \(h_N=o_p(\sqrt{\log n})\), which allows any fixed positive bandwidth. Figure~\ref{fig:tree-localization}(\subref{fig:tree-growing-regions}) shows why \eqref{eq:tree-ancestor-shrinkage} can hold even as we move more levels back from the leaf. In the first tree, moving back one level from depth 3 leaves two shared splits. In the second, moving back two levels from depth 6 leaves four. The shaded intervals show how these additional shared splits restrict the contexts below the ancestor. In general, the number of retained splits, \((d_n-\nu_n)_+\), diverges in probability because \(d_n\) grows as \(\log n\) while \(\nu_n=o_p(\log n)\). Since the tree cycles through the features, each feature range is narrowed repeatedly, by a factor of at most \(1-\theta\) at each split on that feature. Every context below the ancestor therefore lies in a region containing \(x\) whose diameter tends to zero in probability. Continuity gives \(\Lambda_{N,n}(x)\xrightarrow{p}0\).

Our experiments compare leaf and routing weights on the same standard CART partitions. In those trees and forests, the fitting data determine both the splitting feature and its threshold without the spatial balance restriction used above. Verifying \eqref{eq:tree-leaf-count} and \eqref{eq:tree-ancestor-shrinkage} for these fully adaptive CART splits is left for future work.

\paragraph{Forest weights.} The consistency result extends to a forest with a fixed number of trees when each tree satisfies the conditions of Proposition~\ref{prop:routing_consistency}. When the same scenarios are used across the trees, the forest objective is a convex combination of the tree-level objectives. Denoting the objective corresponding to tree \(b\) by \(\widehat F^{(b)}_{N,n}\), we have
\[
\sup_{z\in\mathcal Z}\left|\widehat F^{\mathrm{forest}}_{N,n}(z\mid x)-F(z\mid x)\right|
\leq
\max_{1\leq b\leq n_{\mathrm{trees}}}\sup_{z\in\mathcal Z}\left|\widehat F^{(b)}_{N,n}(z\mid x)-F(z\mid x)\right|
\xrightarrow{p}0.
\]
The consistency of the forest prescription then follows from the same argument as in Corollary~\ref{cor:kernel_consistency}.

\section{Synthetic problem parameters and data-generating processes}
\label{app:synthetic_dgps}

\paragraph{Capacity-expansion parameters.}
For the problem in~\eqref{eq:capacity_expansion_problem}, each product has a distinct primary machine and a secondary machine assigned cyclically. The expansion cost is \(c^E=3\), the production costs are \(c^P_{\mathrm{pri}}=2\) and \(c^P_{\mathrm{sec}}=12\), and the shortage cost is \(c^S=40\). The corresponding production rates are \(\gamma_{\mathrm{pri}}=1\) and \(\gamma_{\mathrm{sec}}=0.7\). Each machine has an initial capacity of \(t_m^0=24\) and a maximum capacity of \(t_m^{\max}=48\).

\paragraph{Feature vectors and demand signals.}
The four demand factors \(U\) and the noise variables \(V\) are independently drawn from \(\mathcal N(0,I_4)\) and \(\mathcal N(0,I_{d_v})\), respectively. All Gaussian draws are then clipped to the interval \([-6,6]\). The core DGPs use the fixed matrix
\[
A=
\begin{pmatrix}
-0.1112 & -0.3822 & 0.4483 & -0.8004\\
0.1177 & -0.7658 & -0.2162 & -0.5941\\
0.0398 & -0.2672 & -0.2241 & -0.9364\\
-0.6477 & 0.7578 & -0.0754 & -0.0218
\end{pmatrix}.
\]
The linear DGP uses \(X=U\) and \(g(U)=AU\). The linear latent DGP uses \(d_v=26\) and
\[
X=Q_{30}^{\top}
\begin{pmatrix}
DU\\
0.25V
\end{pmatrix},
\qquad
D=\operatorname{diag}(2,1.8,1.6,1.4),
\qquad
g(U)=AU.
\]
The nonlinear latent DGPs use \(d_v=12\) and
\[
X=Q_{16}^{\top}
\begin{pmatrix}
U\\
V
\end{pmatrix},
\qquad
g(U)=A\psi(U),
\]
where
\[
\psi(U)=
\begin{pmatrix}
2\sin(2U_1)+U_1U_2\\
1.25(U_2^2-1)+U_2U_3\\
2\sin(2U_3)+U_3U_4\\
1.25(U_4^2-1)+U_4U_1
\end{pmatrix}.
\]
Here, \(Q_{d_x}\) is a fixed orthogonal matrix. In the experiment varying the feature dimension, we change \(d_v\) and the corresponding matrix \(Q_{d_x}\), keeping the four demand factors and the demand signal unchanged. Each component of \(g(U)\) is standardized using the mean and standard deviation estimated from an independent calibration sample of 5000 observations. These calibration values remain fixed across replications.

\paragraph{Demand noise.}
Demand follows \eqref{eq:synthetic_demand_generation}. To specify the noise, define
\[
a(U)=\frac14\sum_{j=1}^4\widetilde g_j(U),
\qquad
b(U;\beta)=0.5+\frac{\beta-0.5}{1+\exp\{-5a(U)\}}.
\]
For \(\varepsilon_j\sim\mathcal N(0,1)\) and \(B_j\sim\operatorname{Bernoulli}(\pi)\), the noise is
\[
\eta_j(U)=
\begin{cases}
2.4\varepsilon_j,
& \text{Homoskedastic},\\[2mm]
3\exp\!\left\{0.75\operatorname{clip}(a(U),-1,1)\right\}\varepsilon_j,
& \text{Conditional variance},\\[2mm]
0.7\varepsilon_j+b(U;40)(B_j-0.15),
& \text{Positive shocks},\\[2mm]
0.7\varepsilon_j+b(U;30)(B_j-0.5),
& \text{Bimodal mixture}.
\end{cases}
\]
The function \(\operatorname{clip}(a,-1,1)\) restricts \(a\) to the interval \([-1,1]\). The positive-shock and bimodal variants use \(\pi=0.15\) and \(\pi=0.5\), respectively, and are used only in Section~\ref{sec:exp_transf}. The Gaussian noise terms and the Bernoulli variables are independent across products and independent of \(U\) and \(V\). In each case, \(\mathbb E[\eta_j(U)\mid U]=0\).

\paragraph{Additional demand information and representation dimension.}
\label{app:additional_information_dgp}
For this experiment, the demand signal is
\[
g(U)=A_{\mathrm{ord}}
\operatorname{diag}\!\left(2,\frac{19}{12},\frac76,\frac34\right)U,
\]
where \(A_{\mathrm{ord}}\) is a separate fixed Gaussian mixing matrix with normalized rows. Demand follows the homoskedastic specification. The noise variables are \(V\sim\mathcal N(0,I_{12})\), independent of \(U\). For the oracle constructions, the first \(k\) demand factors are available at information level \(k=1,\ldots,4\). The learned methods use the full observed feature vector, with the representation dimensions fixed as described in Section~\ref{sec:exp_add_info_and_rep_dimension}.

\section{Training and implementation details}
\label{app:training_details}

\subsection{Synthetic experiments}
\label{app:synt_training_details}

\paragraph{Preprocessing and model fitting.}
The training and validation samples follow Section~\ref{sec:computational_common}. The selected models are not refitted on the combined sample. For each test feature vector, the same conditional demand draws are used to evaluate all methods. \emph{Raw}, \emph{PCA}, \emph{KSPCA}, \emph{MLKR}, and \emph{Neural} standardize the original features using the means and standard deviations from the fitting sample. \emph{KSPCA}, \emph{MLKR}, and \emph{Neural} also standardize demand during fitting.

\emph{PCA}, decision trees, and random forests are implemented using scikit-learn. \emph{KSPCA} uses the degree-two polynomial representation kernel \(\kappa(u,v)=(1+4u^\top v/d_x)^2\) and the linear response kernel \(\ell(\xi,\xi')=\xi^\top\xi'\). \emph{MLKR} follows \citet{mlkr_weinberger}, with an implementation based on the \texttt{metric-learn} package \citep{metriclearn_devazelhes}. We use L-BFGS-B optimization, PCA initialization when feasible, and at most 2000 iterations. The fitted representation scales are retained for \emph{PCA}, \emph{KSPCA}, and \emph{MLKR} when computing the Euclidean distances.

\emph{Neural} is implemented using PyTorch, with ReLU layers of 64, 32, and \(d_\phi\) units, followed by a linear prediction head. It is trained using full-batch AdamW for at most 2000 epochs, with an early-stopping patience of 50. Its learned coordinates are standardized using the means and standard deviations from the fitting sample before the distances are computed. The decision trees and random forests use the multi-output squared-error splitting criterion. Each random forest contains 200 bootstrap trees and considers half of the features at each split. For the leaf weights, if the target leaf contains no historical observations, we move towards the root until reaching an ancestor whose region contains historical observations and assign equal weights to those observations.
\paragraph{Tuning.}
The neural learning rate and weight decay and the tree minimum leaf size are selected using the validation RMSE. Unless fixed by the experiment, the feasible candidate representation dimensions are drawn from \(\{2,\ldots,6\}\) for \emph{PCA} and \emph{MLKR}, \(\{2,3,4\}\) for \emph{KSPCA}, and \(\{2,\ldots,6,16\}\) for \emph{Neural}. For the Gaussian weight constructions, the representation dimension and bandwidth are selected jointly using the validation energy score. The bandwidth candidates multiply a reference distance by
\[
(0.10,0.15,0.20,0.30,0.45,0.65,1.00,1.40,2.00,3.00).
\]

For vector representations, the reference is the median distance from validation feature vectors to their tenth nearest neighbors in the fitting sample; for routing representations, it is the median nearest strictly positive routing distance. Each forest routing distance is divided by twice the maximum depth of its tree. These choices adapt the bandwidth grid to each fitted representation for finite-sample tuning, while Corollary~\ref{cor:kernel_consistency} provides sufficient bandwidth conditions for asymptotic consistency.

Using the \(N_{\mathrm{tr}}\) training observations as the historical scenarios, we compute the energy score for each validation observation \((x,\xi)\) as
\begin{equation}
\mathrm{ES}(x,\xi)=
\sum_{i=1}^{N_{\mathrm{tr}}}w_i(x)\left\|S^{-1}\left(\widetilde{\xi}_i(x)-\xi\right)\right\|_2
-\frac{1}{2}\sum_{i=1}^{N_{\mathrm{tr}}}\sum_{j=1}^{N_{\mathrm{tr}}}w_i(x)w_j(x)\left\|S^{-1}\left(\widetilde{\xi}_i(x)-\widetilde{\xi}_j(x)\right)\right\|_2.
\end{equation}
Here, \(S\) is a diagonal matrix containing the demand standard deviations estimated from the training sample. We select the candidate with the lowest average score across the validation observations \citep{energyscore_gneiting}.

\paragraph{Recentered scenarios.}
LightGBM estimates the conditional mean separately for each demand component, using the validation observations for early stopping. The oracle mean is approximated as described in Section~\ref{sec:exp_transf}. The recentered scenarios are truncated at zero. The representations are learned from the original demand observations, and the representation dimension and bandwidth are selected separately for each scenario construction.

\paragraph{Optimization.} The optimization problems are solved using Gurobi.

\subsection{Kernel comparison}
\label{app:kernel_comparison}

\begin{figure}[H]
    \centering
    \includegraphics[width=0.65\textwidth]{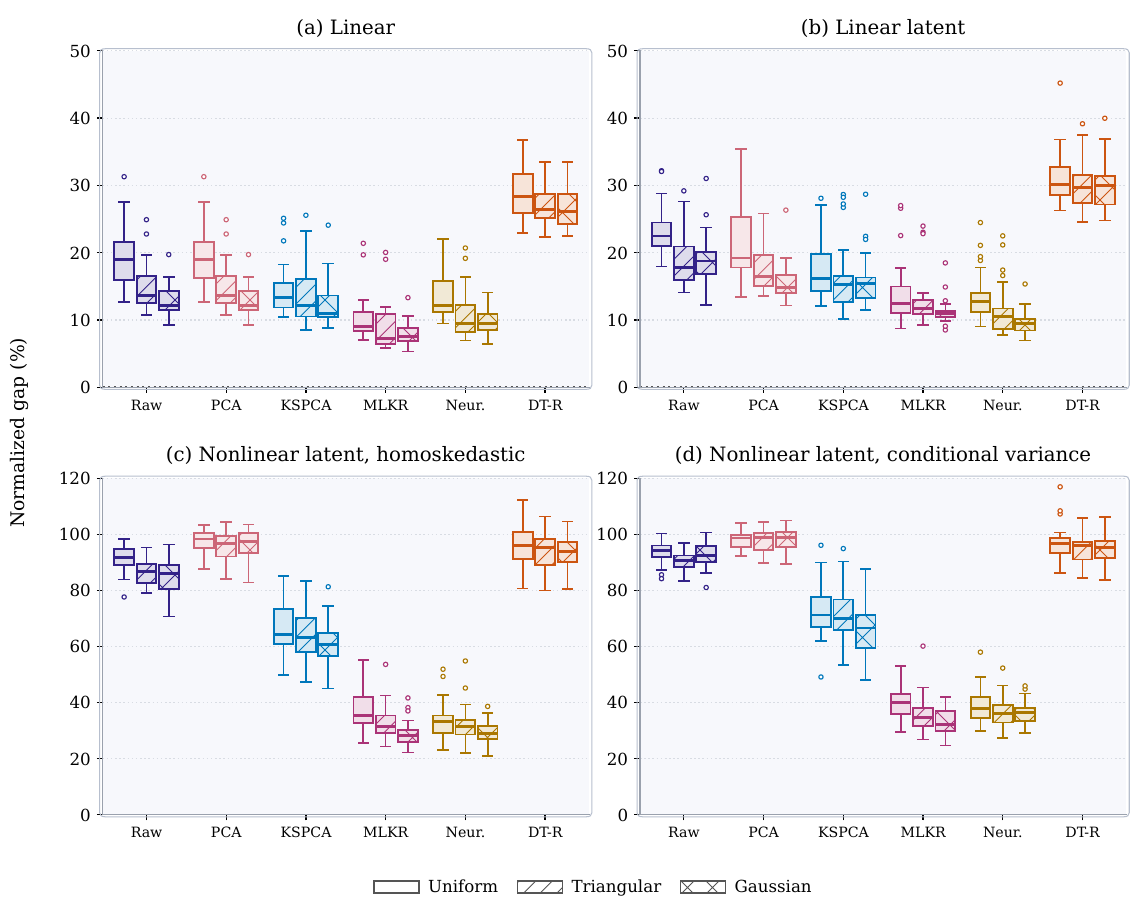}
    \caption{Effect of kernel choice on the normalized gap across the four core DGPs.}
    \label{fig:kernel_choice}
\end{figure}

Figure~\ref{fig:kernel_choice} shows that the Gaussian kernel generally achieves slightly lower normalized gaps than the uniform and triangular kernels, although the differences are modest and vary across methods. This may reflect the ability of the Gaussian kernel to distinguish among nearby observations, unlike the uniform kernel, while assigning smoothly decreasing weights to more distant observations without imposing a hard cutoff. Overall, the choice of feature representation has a substantially larger effect on prescription quality than the choice of kernel.

\subsection{Comparison at a common representation dimension}
\label{app:common_dimension}

\begin{figure}[H]
    \centering
    \includegraphics[width=0.65\textwidth]{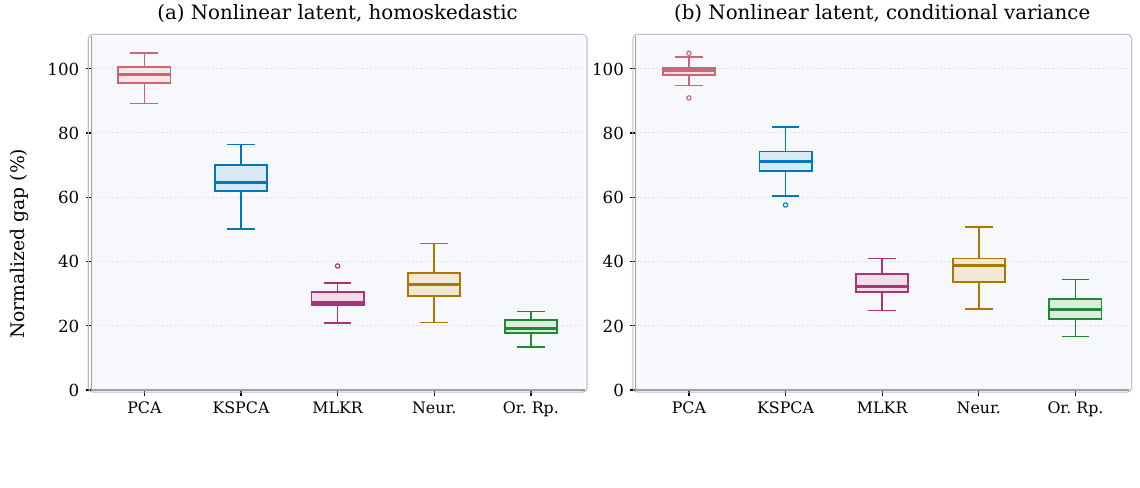}
    \caption{Normalized gap in the two nonlinear DGPs with the representation dimension fixed at \(d_\phi=4\).}
    \label{fig:exp2b}
\end{figure}

Figure~\ref{fig:exp2b} compares the methods at a common representation dimension of \(d_\phi=4\) in the two nonlinear DGPs. Substantial differences in prescription quality remain, indicating that the choice of representation-learning method matters beyond the representation dimension alone.

\subsection{Larger learned representations}
\label{app:larger_representations}

\begin{figure}[H]
    \centering
    \includegraphics[width=0.65\textwidth]{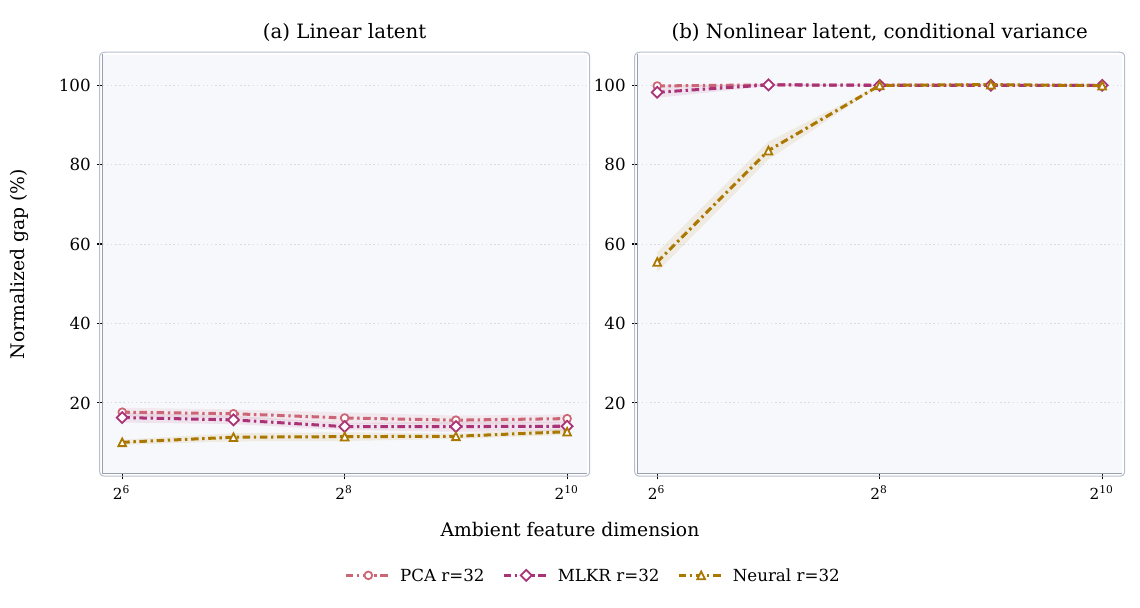}
    \caption{Effect of the feature dimension on \emph{PCA}, \emph{MLKR}, and \emph{Neural} with \(d_\phi=32\).}
    \label{fig:exp4b_r32}
\end{figure}

Figure~\ref{fig:exp4b_r32} shows the results for \emph{PCA}, \emph{MLKR}, and \emph{Neural} with a representation dimension of \(d_\phi=32\). In the linear latent setting, all three methods maintain relatively low gaps as the feature dimension increases. In the nonlinear setting, \emph{PCA} and \emph{MLKR} perform close to \emph{SAA}, while the larger \emph{Neural} representation avoids performing worse than \emph{SAA} at the largest feature dimensions, without recovering an improvement over it.

\section{Details on the contextual electricity dispatch problem}
\label{app:gefcom_details}

\subsection{Two-stage model and parameters}

We use the two-stage formulation adapted from \citet{zakeri2019pricing}. There are no constraints linking the hours, so the model below describes one fixed hour. The index \(t\) denotes the day.

Let \(\mathcal P\), \(\mathcal B\), and \(\mathcal E\) denote the generators, network nodes, and transmission lines. For node \(k\), let \(\mathcal P_k\) contain its generators, and let \(\mathcal E_k^{\mathrm{in}}\) and \(\mathcal E_k^{\mathrm{out}}\) contain the incoming and outgoing lines under a fixed orientation. With generation targets \(z_g\) and capacities \(\bar p_g\), the contextual problem is
\[
\min_{\substack{0\leq z_g\leq\bar p_g\\g\in\mathcal P}}
\sum_{i=t-365}^{t-1}w_i(X_t)
Q^{\mathrm{disp}}\!\left(z,\widetilde\xi_i^{\mathrm{hr}}(X_t)\right),
\]
where \(\widetilde\xi_i^{\mathrm{hr}}(X_t)\) is the component of the scenario from day \(i\) corresponding to the fixed hour.

For the realized hourly demand \(\xi\geq0\), let \(p_g\) denote actual generation, \(p_g^+\) and \(p_g^-\) the upward and downward adjustments, \(p_e^{\mathrm{line}}\) the signed line flow, and \(p^{\mathrm{shed}}\) the unmet system demand. With nominal node demands \(\xi_k^{\mathrm{base}}\), the recourse problem is
\[
\begin{aligned}
Q^{\mathrm{disp}}(z,\xi)=\min\quad
&\sum_{g\in\mathcal P}
\left[
c_g^{(2)}p_g^2+c_g^{(1)}p_g
+c_g^+p_g^++c_g^-p_g^-
\right]
+c^{\mathrm{shed}}p^{\mathrm{shed}}\\
\text{s.t.}\quad
&p_g=z_g+p_g^+-p_g^-,
&&g\in\mathcal P,\\
&\sum_{g\in\mathcal P_k}p_g
+\sum_{e\in\mathcal E_k^{\mathrm{in}}}p_e^{\mathrm{line}}
-\sum_{e\in\mathcal E_k^{\mathrm{out}}}p_e^{\mathrm{line}}
=\frac{\xi_k^{\mathrm{base}}}
{\sum_{j\in\mathcal B}\xi_j^{\mathrm{base}}}
\left(\xi-p^{\mathrm{shed}}\right),
&&k\in\mathcal B,\\
&0\leq p_g\leq\bar p_g,\qquad
p_g^+\geq0,\quad p_g^-\geq0,
&&g\in\mathcal P,\\
&p^{\mathrm{shed}}\geq0,\qquad
p^{\mathrm{line}}\in\mathcal F.
\end{aligned}
\]
The minimization is over the generation, adjustment, line flow, and unmet-demand variables. The set \(\mathcal F\) imposes the lossless DC network constraints and the transmission limits \(-\bar p_e^{\mathrm{line}}\leq p_e^{\mathrm{line}}\leq\bar p_e^{\mathrm{line}}\). The node balances allocate demand and unmet demand according to the nominal demand shares. The coefficients \(c_g^{(2)}\) and \(c_g^{(1)}\) specify the production costs, \(c_g^+\) and \(c_g^-\) the adjustment costs, and \(c^{\mathrm{shed}}\) the penalty for unmet demand.

The network has 30 nodes, six generators, and 41 lines from MATPOWER's \texttt{case30} \citep{zimmerman2011matpower}. The generator capacities and line limits are multiplied by the smallest positive integer that accommodates the maximum system demand observed during 2004--2007 under the DC network constraints. The quadratic production-cost coefficients are divided by the same factor, while the linear coefficients remain unchanged. We set
\[
c_g^+=5\left(c_g^{(2)}\bar p_g+c_g^{(1)}\right),
\qquad
c_g^-=\tfrac14c_g^+,
\qquad
c^{\mathrm{shed}}=500.
\]
The daily cost is the sum of the 24 hourly costs. The optimization problems are solved using Gurobi.

\subsection{Data construction}
\label{app:dispatch_data}

The GEFCom2012 data contain hourly electricity demand for 20 zones and temperatures from 11 weather stations \citep{hong2014gefcom}. We restore the withheld demand observations using the published competition solution and aggregate the zonal demand into a vector \(\xi_t\in\mathbb R^{24}\) of hourly system demand for each day \(t\). Each missing temperature is replaced by the most recent available observation for the same station and hour.

The feature vector \(X_t\) contains the preceding seven days of hourly zonal and system demand and temperatures, together with the temperatures for day \(t\). It also includes weekday and month indicators, a holiday indicator, and a linear time trend. For each hour, additional features include demand lags of one, two, and seven days, the corresponding seven-day mean and standard deviation, and the mean, minimum, and maximum temperature across the weather stations. The previous day's mean, minimum, maximum, and standard deviation of system demand provide four further features.

The feature vector contains 6004 components, comprising 5640 demand and temperature features, 168 calendar features, 192 additional hourly features, and four daily demand statistics. The realized temperatures for day \(t\) are used retrospectively as proxies for the day-ahead temperature forecasts. After accounting for the seven-day history, the dataset contains 1635 daily observations from 8 January 2004 to 29 June 2008. Of these, 1089 observations from 2004--2006 are used for initial fitting, 365 from 2007 for validation, and 181 from 2008 for evaluation.

\subsection{Representation fitting and scenario construction}
\label{app:dispatch_training}

After tuning on the 2007 validation data, the selected models are refitted using the observations from 2004--2007. Unless stated otherwise, the preprocessing and random forest settings follow Appendix~\ref{app:synt_training_details}.

\emph{MLKR} and \emph{Neural} consider \(d_\phi\in\{4,8,16,24,32\}\), while \emph{KSPCA} considers the same dimensions up to 24 and uses the polynomial representation kernel specified above. For \emph{MLKR}, PCA first reduces the 5640 demand and temperature features to 256 components, which are then combined with the remaining 364 features. The training allows at most 20,000 iterations. \emph{Neural} uses ReLU layers with 128, 64, and \(d_\phi\) units and full-batch AdamW training for at most 300 epochs, with an early-stopping patience of 25.

\emph{TSMixer} \citep{chen2023tsmixer} uses the \texttt{pytorch-tsmixer} implementation from the repository (\url{https://github.com/ditschuk/pytorch-tsmixer}) with three mixer blocks, 32 hidden channels, a feed-forward width of 64, and a dropout of 0.1. Training uses AdamW with batches of 128 for at most 150 epochs, with an early-stopping patience of 20. Both neural architectures minimize the squared prediction error on the standardized demand. The learning rate and weight decay are selected using the validation RMSE, and the final training uses the number of epochs selected during validation. The \emph{TSMixer} representation contains the \(24\times32\) hidden activations before the prediction head, flattened into 768 coordinates and standardized using the means and standard deviations from the fitting sample. TSMixer uses the same contextual information as the other methods, organizing the historical observations and the target-day features along the temporal and feature dimensions of the architecture.

LightGBM fits one model per hour to estimate the conditional mean from the complete feature vector \citep{ke2017lightgbm}, using a squared-error objective and the validation observations for early stopping.

For LightGBM and TSMixer recentering, historical predictions are generated using a rolling-origin procedure in successive six-month blocks. For each block, the preceding six months are used for validation, after which the selected model is refitted using all the observations preceding the block. If \(\widehat{\xi}_t\) denotes the vector of predicted hourly demands for day \(t\), the scenario from day \(i\) is
\[
\widetilde\xi_i(X_t)
=
[\widehat\xi_t+\xi_i-\widehat\xi_i]_+,
\]
where the positive-part operator is applied componentwise. The representations are learned from the original demand observations. The Gaussian weight constructions are tuned separately for each scenario construction using the validation energy score over rolling windows of 365 days in 2007. The candidate representation dimensions are used where applicable, together with the bandwidth multipliers in Appendix~\ref{app:synt_training_details}. The reference distance is the median distance to the tenth nearest neighbor across the validation days. After the final representation is fitted, this distance is recomputed and the selected multiplier is retained. The evaluation then follows the rolling-window procedure in Section~\ref{sec:electricity_dispatch}.

\end{document}